\documentclass[12pt,dvips]{amsart}
\usepackage{euler, amsfonts, amssymb, latexsym, epsfig,epic,wasysym}

\newcommand\lie[1]{{\mathfrak #1}}

\newcommand\Tr{{\rm Tr\,}}
\newcommand\diag{{\rm diag}}

\newcommand\iso{{\ \cong\ }}
\newcommand\tensor{{\otimes}}

\newtheorem{Theorem}{Theorem} 
\newtheorem{Proposition}{Proposition} 
\newtheorem{Lemma}{Lemma}

\newtheorem{Corollary}{Corollary}
\newtheorem*{Corollary*}{Corollary}
 
\newtheorem*{Theorem*}{Theorem}
\theoremstyle{remark}

\newcommand\onto{\mathop{\twoheadrightarrow}}
\newcommand\into{\operatorname*{\hookrightarrow}}

\newcommand\Pone{{\mathbb P}^1}

\newcommand\complexes{{\mathbb C}}
\newcommand\integers{{\mathbb Z}}

\newcommand\onehalf{\frac{1}{2}}
\newcommand\naturals{{\mathbb N}}

\newcommand\<{\langle}
\renewcommand\>{\rangle}
\newcommand\junk[1]{}
\newcommand\To{\longrightarrow}

\theoremstyle{plain}

\newcommand\dfn{\bf} % maybe should be \em

\newcommand\seclabel[1]{\label{sec:#1}}

\newcommand\cycinside[1]{#1}
\newcommand\cyc[1]{\circlearrowleft (\cycinside{#1})}
\newcommand\cp{\bullet} % new product
\newcommand\der{\partial}
\newcommand\mdeg{{\rm m}\!\deg}
\newcommand\codim{{\rm co}\!\dim}
\newcommand\SN{S_N}
\newcommand\MNC{{M_N(\complexes)}}

\newcommand\MNCt{{M_N(\complexes[t])}}
\newcommand\A{A}

\begin{document}
\pagestyle{plain}

\title{A scheme related to the Brauer loop model}
\author{Allen Knutson}
\address{Allen Knutson, 970 Evans Hall, 
  University of California at Berkeley, Berkeley CA 94720}
\email{allenk@math.berkeley.edu}
\thanks{AK was supported by NSF grant 0303523.}
\author{Paul Zinn-Justin}
\address{Paul Zinn-Justin, Poncelet Laboratory (UMI 2615 of CNRS), 
  Independent University of
  Moscow and LPTMS (UMR 8626 of CNRS), Universit\'e Paris-Sud.}
\email{pzinn@lptms.u-psud.fr}
\thanks{PZJ was supported by the ENIGMA MRTN-CT-2004-5652 European network 
  and MISGAM ESF program}
\date{\today}

\begin{abstract}
  We introduce the {\dfn Brauer loop scheme} $E := \{M\in \MNC : M\cp M = 0\}$,
  where $\cp$ is a certain degeneration of the ordinary matrix product.
  Its components of top dimension, $\lfloor N^2/2\rfloor$, correspond
  to involutions $\pi \in \SN$ having one or no fixed points.
  In the case $N$ even, this scheme contains the upper-upper scheme 
  from [Knutson '04] as a union of $(N/2)!$ of its components.
  One of those is a degeneration of the {\em commuting variety} of
  pairs of commuting matrices.
  
  The {\em Brauer loop model} is an
  integrable stochastic process studied in [de Gier--Nienhuis '04],
  based on earlier related work in [Martins--Nienhuis--Rietman '98],
  and some of the entries of its
  Perron--Frobenius eigenvector were observed (conjecturally) to equal 
  the degrees of the components of the upper-upper scheme.

  Our proof of this equality follows the program outlined in
  [Di Francesco--Zinn-Justin '04]. In that paper, the entries of the
  Perron--Frobenius eigenvector were generalized from numbers to
  polynomials, which allowed them to be calculated inductively using
  divided difference operators.  We relate these polynomials to the
  {\em multi}degrees of the components of the Brauer loop scheme,
  defined using an evident torus action on $E$. As a consequence, 
  we obtain a formula for the degree of the commuting variety,
  previously calculated up to $4\times 4$ matrices.
\end{abstract}

\maketitle

{ \small\tableofcontents}
%
%%%%%%%%%%%%%%%%%%%%%%%%%%%%%%%%%%%%%%%%%%%%%%%%%%%%%%%%%%%%%%%%%%%%%
%

\section{Introduction}\seclabel{intro}

\subsection{The scheme $E$}

Let $N$ be a positive integer. 
Call a sequence $(i_1,\ldots,i_k) \in \{1,\ldots,N\}^k$
{\dfn cyclically ordered}, written ``$\cyc{i_1\le\ldots\le i_k}$'',
if for some $j$, the rotated sequence
$(i_j, i_{j+1}, \ldots, i_k, i_1, \ldots, i_{j-1})$ is weakly increasing,
with the additional condition that all $i$'s are equal if $i_1=i_k$. 

We define a new product $\cp$ on the space $\MNC$ of $N\times N$
matrices according to the rule
\begin{equation*}
  (P\cp Q)_{ik}=\sum_{j:\ \cyc{i\le j\le k}} P_{ij} Q_{jk}\qquad i,k=1,\ldots,N
\end{equation*}
This is a degeneration of the usual product, as we explain in section
\ref{ssec:algdegen}. Hence it is associative, and indeed has triple product
\begin{equation*}\label{eqn:triple}
  (P\cp Q\cp R)_{il}=\sum_{j,k:\ \cyc{i\le j\le k\le l}} P_{ij} Q_{jk} R_{kl}.
\end{equation*}
With the usual addition, $(\MNC,\cp)$ forms an algebra,
and the identity matrix is again the unit. 
A matrix $P$ possesses an % (left, right)
inverse $P^{\cp-1}$ for this product
if and only if its diagonal entries are all non-zero. In particular,
the set
$$ U = \{M \in \MNC : M_{ii} = 1, i=1\ldots N \} $$
forms a group under $\cp$. We note that if $R$, $R'$ are upper triangular,
then $R\cp R'=RR'$ and $R^{\cp-1} = R^{-1}$.
All these facts are easiest to see in the model of $\cp$ presented
in section \ref{ssec:semidirect}.

The {\dfn cycling automorphism} $M'_{ij} := M_{i+1,j+1}$, where the indices
are taken mod $N$, is an automorphism of both the ordinary and $\cp$
multiplications. This automorphism is inner for the ordinary multiplication,
but is an outer automorphism for $\cp$. With this automorphism in mind,
essentially every reference to $i$, $j$, etc. in this paper has an implicit
``$\bmod\ N$''.

We define the {\dfn Brauer loop scheme} $E$ to be the space of matrices 
$M\in \MNC$ that satisfy $M\cp M=0$, and have zero diagonal.\footnote{%
  These equations are not redundant: $M\cp M=0$ implies that $M_{ii}^2=0$
  for each $i$, but not that $M_{ii}=0$. This is an empty distinction on the 
  set but an important distinction on the scheme, and one that affects the
  (multi)degree that will interest us later. A similar phenomenon occurs
  already with ordinary matrix multiplication
  (as we address in section~\ref{sec:flatlimit}): while any matrix with
  $M^2=0$ has zero trace, the linear trace condition can't be inferred
  algebraically from the quadratic conditions $M^2=0$.}
(The name will be explained in section \ref{ssec:brauerpolys}.)  
In equations, we require
\begin{eqnarray*}
  \sum_{\cyc{i\le j\le k}} M_{ij} M_{jk} 
  &=& 0 \qquad i,k=1,\ldots,N, i\neq k \\
  M_{ii} &=& 0 \qquad i=1,\ldots,N 
\end{eqnarray*}

The scheme $E$ looks similar to the irreducible scheme $\{M : M^2=0\}$
(a precise relation is spelled out in section \ref{sec:flatlimit}),
and in particular has the same dimension $\lfloor N^2/2\rfloor$.
However, $E$ is reducible, and we now describe its components of
top dimension.

In what follows the parity of $N$ will play a role, so write 
$$ N=2n+r, \qquad  r=0 \mbox{ or } 1.$$
We will refer to involutions of $\{1,\ldots,N\}$ with $r$ fixed
points as {\dfn link patterns}, and draw them as chord diagrams in the disk. 
In particular, the $2$-cycles of an involution will be referred to as chords,
and a ``crossing'' in a link pattern is a pair of chords which
cross each other when drawn as segments in the disk.
There are $(N-1)!! := (N-1)(N-3)(N-5)\cdots(1+r)$ link patterns of size $N$.

\noindent {\it Example}: The involution with cycles $(1 5)$ $(2 4)$ $(3 6)$ is represented
as \begin{minipage}{1.8cm}\hfil\epsfig{file=diag3-11.eps, width=1.5cm}\end{minipage}.

The following is a combination of theorems \ref{thm:diagpaired} and
\ref{thm:topcomps}.

\begin{Theorem*} 
  For each $M \in E$, the nonzero elements of the diagonal of $M^2$ 
  (with respect to {\em ordinary} multiplication) come in equal pairs.
  Put another way, there is a link pattern $\pi$
  such that $(M^2)_{ii} = (M^2)_{\pi(i)\pi(i)}$ for all $i$.
  In addition, $(M^2)_{ii} = 0$ if $\pi(i)=i$.

  Conversely, for each such $\pi$, the open subscheme
  $$ \big\{ M\in E 
  : (M^2)_{ii} = (M^2)_{jj} \hbox{ if and only if } j\in\{i,\pi(i)\} \big\} $$
  is nonempty, irreducible, 
  and of dimension $\lfloor N^2/2 \rfloor$.
\end{Theorem*}

Hence, each $E_\pi$ defined by
$$ E_\pi := \overline{
   \big\{ M\in E 
  : (M^2)_{ii} = (M^2)_{jj} \hbox{ if and only if } j\in\{i,\pi(i)\} \big\}}$$
is a component of $E$.
In fact we conjecture that $E = \cup_\pi E_\pi$.\footnote{This has now
  been proven by Brian Rothbach; details will appear elsewhere.}
The closest we come to proving this, in theorems \ref{thm:topcomps} and
\ref{thm:reduced}, is

\begin{Theorem*}
  If $E \neq \cup_\pi E_\pi$, where $\pi$ runs over the set of
  link patterns, then $\dim (E \setminus \cup_\pi E_\pi) < \dim E$.
  Also, $E$ is generically reduced along each $E_\pi$.
\end{Theorem*}

Theorem \ref{thm:reduced} gives a different characterization of
the $\{E_\pi\}$:

\newcommand\Mpi{\underline \pi}
\newcommand\Mrho{\underline \rho}
\newcommand\Mpisl{\Mpi_<}
\begin{Theorem*}
  Let $\Mpi$ denote the permutation matrix of a link pattern $\pi$,
  with the diagonal zeroed out if $\pi$ has a fixed point (i.e. if $N$
  is odd). Then
  $$ E_\pi = \overline{U \cdot \{t \Mpi : t \hbox{ diagonal} \} } $$
  where $U$ acts by $\cp$-conjugation.
\end{Theorem*}

This lets us determine in theorem \ref{thm:compeqns} some (and
conjecturally, all) of the defining equations of the $\{E_\pi\}$.

The cycling automorphism acts on $E$, and on the set of link
patterns by rotation.  
We will make use, too, of the action of the full symmetric group $\SN$
on the set % $CP_N$ 
of link patterns by conjugation, even though $\SN$ does not act on $E$.
For $i=1,\ldots,N$,
denote by $f_i$ the transposition $i\leftrightarrow i+1$ (where $N+1\equiv 1$),
and let $f_i\cdot\pi := f_i\circ\pi\circ f_i^{-1}$.

For each $i=1\ldots N$, there is an idempotent 
{\dfn Temperley--Lieb operator} $e_i$ on the
set of link patterns defined by
$$ (e_i\cdot \pi)(j) =
\begin{cases}
  i+(i+1)-j &\mbox{ if $j=i$ or $j=i+1$ } \\
  \pi(i)+\pi(i+1)-j &\mbox{ if $j=\pi(i)$ or $j=\pi(i+1)$ } \\
  \pi(j) &\mbox{ otherwise }  
\end{cases}
$$
where all addition is mod $N$. Graphically, $e_i$ connects the chords
coming to $i,i+1$ to one another, and puts in a new chord connecting $i,i+1$.

Together, the $\{f_i\}$ and $\{e_i\}$ form a representation of the
affine {\dfn $O(1)$ Brauer algebra}.  (Actually, the $f_i$ and $e_i$,
$i=1,\ldots,N-1$, which satisfy the relations of the usual $O(1)$
Brauer algebra, are enough to generate the whole algebra of
operators.)  The Brauer algebra is itself a degenerate point of the
braid-monoid algebra.

\subsection{The Brauer loop model polynomials $\{\Psi_\pi\}$}
\label{ssec:brauerpolys}

In \cite{dGN} there is associated to each link pattern $\pi$ a positive integer
$d_\pi$, as follows. (They will at first only appear to be rationals.)

\begin{figure}[htbp]
  \centering
  \epsfig{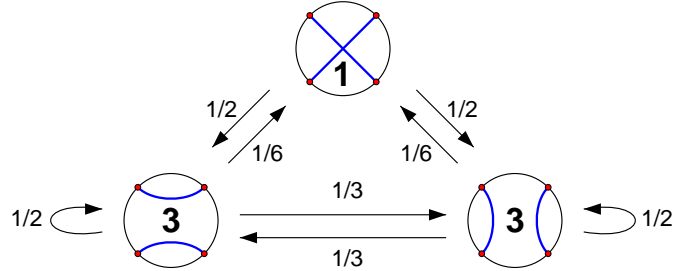}
  \caption{The de Gier--Nienhuis Markov process for $N=4$. 
    The edges are labeled with transition probabilities, and the 
    nodes with the stationary distribution, rescaled to make the 
    minimum value $1$.}
  \label{fig:markov4}
\end{figure}

Consider a Markov process whose states are the set of link patterns.
The transitions from a link pattern $\pi$ are
to $\{e_i \cdot \pi, f_i \cdot \pi\}$,
where $i$ is chosen with equal probability from $1\ldots N$, 
and $e_i,f_i$ are then chosen with probabilities $2/3$ and $1/3$ 
(see figure~\ref{fig:markov4}). For the origin of this Markov process and
its relation to standard quantum integrable models, see \cite{MR,MNR}.

Many conjectures are stated in \cite{dGN}, among them that
the least probable link patterns are those with the most crossings,
and all other probabilities are {\em integer} multiples $d_\pi$ 
of the least probable. Some of these $d_\pi$ were also noticed to
match the degrees of the components of the upper-upper scheme from
\cite{Kn}. (That scheme reappears here in section \ref{sec:perm}.)

This integrality, and much else, was proven in \cite{DFZJ} by
considering a slightly different Markov process, and generalizing
the $d_\pi$ to polynomials $\Psi_\pi \in \integers[z_1,\ldots,z_N]$.
(In fact \cite{DFZJ} only treats the case $N$ even, but as indicated
in \cite{dGN} the case $N$ odd is very similar.)  

In this more general Markov process, the transition amplitudes are given 
by the so-called {\em transfer matrix}, one possible definition of which is
\begin{equation*}
  T'=\prod_{i=1}^{n} \prod_{j=1}^{n} 
  {\check R}_{i+2j-2}(z_{2j-1} - z_{2i+2j-2}), \quad\mbox{ where }
  \check R_i(u) = a(u) 1 + b(u) f_i + c(u) e_i
\end{equation*}
(for $N=2n$; a similar expression exists for $N$ odd)
where the product is ordered and the indices meant mod $N$.  
Here
$\check R_i(u)$ is a linear combination of ${\bf 1}$, $f_i$, and $e_i$ (their 
action on link patterns extended %by linearity to the whole vector space), 
to turn them into linear operators), with coefficients dependent on
$u$ still to be determined (as they will be below, in equation 
(\ref{eqn:Rcheckform})).

The stationary distribution can be encoded into a vector 
$\Psi=\sum_\pi \Psi_\pi\ \pi$ in the space of 
$\integers[z_1,\ldots,z_n]$-linear combinations of link patterns. 
Again, we scale $\Psi$ to clear denominators, making the $\{\Psi_\pi\}$
polynomials with \hbox{$GCD=1$}. The unnormalized 
probabilities $d_\pi$ of the de Gier--Nienhuis
process are recovered by setting all $z$'s to zero.

We now explain the coefficients we will use in $\check R_i$.
Consider the equations
\begin{equation} \label{eqn:Rcheck}
  \check R_i(z_i-z_{i+1}) \Psi = \tau_i \Psi \qquad i=1,\ldots,N
\end{equation}
where $\tau_i$ switches the variables $z_i$ and $z_{i+1}$, i.e.
$\tau_i F(\ldots,z_i,z_{i+1},\ldots)=F(\ldots,z_{i+1},z_i,\ldots)$.
Since $T'$, by rearranging its definition,
is a product of $\tau_i\check R_i(z_i-z_{i+1})$ operators, 
the equations~(\ref{eqn:Rcheck}) imply that $T' \Psi=\Psi$. 
More generally, products of $\tau_i \check R_i$, with the condition that the
corresponding product of $\tau_i$ is the identity (which ensures that these operators
have well-defined meaning once the $z_i$ are specialized to complex numbers),
generate a whole
algebra of operators acting on $\integers[z_1,\ldots,z_N]$-linear 
combinations of link patterns.
The ``integrability'' condition implies
that this algebra must be commutative;
this is achieved by requiring
that the $\check R_i$ generate a representation of the affine
$\hat{A}_{N-1}$ Weyl group, or equivalently that
%The $\pi$-component of this equation reads
%\begin{equation} \label{eqn:picomp}
%a(z_i-z_{i+1}) \Psi_\pi + b(z_i-z_{i+1}) \Psi_{f_i\cdot \pi} 
%+ c(z_i-z_{i+1}) \sum_{\rho:\ e_i\cdot \rho=\pi} \Psi_\rho
%= \tau_i \Psi_\pi.
%\end{equation}
%
%Applying the equations $\tau_i^2=1$ and 
%$\tau_i\tau_{i+1}\tau_i = \tau_{i+1}\tau_i\tau_{i+1}$ to $\Psi$,
%and using equation (\ref{eqn:Rcheck}), we find
%$$ \check R_i(z_{i+1}-z_i) \check R_i(z_i-z_{i+1}) \Psi = \Psi $$
%and
%$$
%\check R_i(z_{i+1}-z_{i+2})\check R_{i+1}(z_i-z_{i+2})\check R_i(z_i-z_{i+1})
%\Psi 
%= \check R_{i+1}(z_i-z_{i+1})\check R_i(z_i-z_{i+2})
%\check R_{i+1}(z_{i+1}-z_{i+2}) \Psi. $$ 
%If we ask that the operators $\tau_i \check R_i(z_{i+1}-z_i)$ be 
%simultaneously diagonalizable, rather than merely have a common fixed vector
%$\Psi$, we must impose the conditions
$ \check R_i(u)\check R_i(-u)=I$ 
and 
$\check R_i(u)\check R_{i+1}(u+v)\check R_i(v)
=\check R_{i+1}(v)\check R_i(u+v)\check R_{i+1}(u)$ 
(the {\em Yang--Baxter equation}).

As shown in \cite{NR}, these two equations on the $\{\check R_i\}$
fix uniquely the coefficients $a(u)$, $b(u)$, $c(u)$ in their definition 
(up to scaling of $u$, and multiplication
by a function $\phi(u)$ satisfying $\phi(u)\phi(-u)=1$), resulting in
\begin{equation}
  \label{eqn:Rcheckform}
   \check R_i(u)=\left( (1-u) I + \frac{u}{2}(1-u) f_i + u e_i\right)
   \bigg /\left((1-\frac{u}{2})(1+u)\right). 
\end{equation}

These polynomials $\{\Psi_\pi\}$ are characterized %up to an overall sign
by the following two properties
(equations (3.11--14) and (3.19) of \cite{DFZJ}): 

\begin{itemize}
\item
  Recall that $N=2n+r$, $r=0,1$, and define $\pi_0(i)=i+n \bmod 2n$ for
  $i\leq 2n$, and $\pi_0(N)=N$ if N is odd. So $\pi_0$ is a
  {\dfn maximally crossing link pattern}, and the only such if $N$ is even. 
  Then
  \begin{equation} \label{eqn:pizeroa}
      \Psi_{\pi_0} 
      =\prod_{i=1\ldots N \atop j:\ \cyc{i<j<i+n}} \!\!\!\!\!\! (1+z_i-z_j)\ \ 
      \left( \prod_{i=n+1}^N (1+z_i-z_{i+n}) \right)^r
  \end{equation}
\item 
  If $\pi(i) \neq i+1$, then 
  % $\Psi_{\pi}/(1+z_i-z_{i+1})$ is a polynomial,  and
  \begin{equation}
    \label{eqn:schub}
    \Psi_{\pi}
    +\Psi_{f_i\cdot\pi}=
    -\frac{2+z_{i+1}-z_{i}}{1+z_{i+1}-z_i} 
    \,\der_i 
    \left((1+z_{i+1}-z_i)\Psi_{\pi}\right)
  \end{equation}
  where $\partial_i F := (F - \tau_i F)/(z_i - z_{i+1})$.
\end{itemize}

The equations (\ref{eqn:schub}) allow one to express a component in terms of
another with one more crossing (and, for $N$ odd, to move the fixed point),
until one reaches $\Psi_{\pi_0}$, which is given by (\ref{eqn:pizeroa}).  
They are simply the components of equation
(\ref{eqn:Rcheck}) for which $\pi(i)\ne i+1$. 
Equation (\ref{eqn:pizeroa}) is also required by the compatibility of
the set of equations (\ref{eqn:Rcheck}) and the coprimality of the
$\Psi_\pi$, see \cite{DFZJ}.

%They are derived by
%The equations (\ref{eqn:schub}) come from 
%plugging equation (\ref{eqn:Rcheckform}) for $a(),b(),c()$ 
%into equation (\ref{eqn:picomp}), where the summation term disappears
%because $\pi(i)\neq \pi(i+1)$.

The purpose of this paper is to give a geometric interpretation of these
polynomials, extending and proving the observations of \cite{dGN} (at
$z_i\equiv 0$) and of \cite{DFZJ}.  Having one allows us to obtain a
positivity property of the $\{\Psi_\pi\}$, stated
below (corollary \ref{cor:positive}).

\subsection{Degrees and multidegrees}

Since $E$ (and each $E_\pi$) is invariant under rescaling, it is the
affine cone over a projective scheme, and thus has a well-defined
degree.  The {\dfn degree} $\deg_W X \in \naturals$ of an affine cone
$X$ in a vector space $W$ is characterized by three properties:
\begin{enumerate}
\item[1.] If $X=W=\{0\}$, then $\deg_W X = 1$.
\item[2.] If the scheme $X$ has top-dimensional components $X_i$, 
  where $m_i>0$ denotes the multiplicity of $X_i$ in $X$, 
  then $\deg_W X = \sum_i m_i \deg_W X_i$. This lets one reduce from
  the case of schemes to the case of varieties (reduced irreducible schemes).
\item[3.] If $X$ is a variety, and $H$ is a hyperplane in $W$, 
  then $\deg_W X = \deg_H (X\cap H)$.
  (Note that $X\cap H$ may be neither reduced nor irreducible.)
\end{enumerate}
From these it is easy to see that the degree is a nonnegative integer
(and only zero if $X$ is empty); more work is necessary to know that
it is well-defined, but of course this is standard.

\begin{Theorem*}
  For each link pattern $\pi$
  the probability of the state $\pi$ in
  the de Gier--Nienhuis Markov process is proportional to the degree
  of the component $E_\pi$. More precisely, $\deg E_\pi=d_\pi$.
\end{Theorem*}

This was conjectured in \cite{dGN} for those $\pi$ satisfying
$i \leq n \implies \pi(i)>n$
(see section \ref{sec:perm} for the meaning of this condition).
Most elements of a proof in that special case were given in \cite{DFZJ}, 
by going beyond degrees to the more refined {\em multidegrees} of the
components $\{E_\pi\}$ of $E$ (our reference for multidegrees is \cite{MS}).
These are defined using an additional torus action on $E$, 
the conjugation action by invertible diagonal matrices,
with respect to either $\cp$ or ordinary multiplication (the action
is the same). Let $T$ denote the $(N+1)$-dimensional torus
$$ T := \hbox{the rescaling $\complexes^\times$} \ \times\ 
\hbox{the invertible diagonal matrices in $\MNC$}. $$
If we denote the usual basis for $T$'s weight lattice by $(\A,z_1,\ldots,z_N)$,
then the weights of the $T$-action on $\MNC$ are 
$\{\A + z_i - z_j\}, i,j = 1\ldots N$.

When $T$ acts on a vector space $W$
preserving a subscheme $X$, one has an associated homogeneous
{\dfn multidegree} $\mdeg_W X \in \integers[A,z_1,\ldots,z_N]$.
It is also known as the equivariant cohomology class, 
and the equivariant Chow class \cite{Br}.\footnote{%
  It is closely related to the equivariant multiplicity \cite{Ro}, 
  which is best thought of homologically rather than cohomologically.
  The equivariant multiplicity was introduced in \cite{Jo} to study 
  {\em orbital varieties}, the components of the intersection of
  a nilpotent orbit with the upper triangular matrices. In this paper
  we degenerate the nilpotent orbit $\{M^2=0\}$ to get components, 
  rather than intersect with anything.}
The multidegree is characterized by similar axioms to the degree, 
except that the third one is split into two cases:

\begin{enumerate}
\item[3'.] Assume $X$ is a variety, 
  and $H$ is a $T$-invariant hyperplane in $W$.
  \begin{enumerate}
  \item If $X\not\subset H$, then $\mdeg_W X = \mdeg_H (X\cap H)$.
  \item If $X\subset H$, then 
  $ \mdeg_W X = (\mdeg_H X) \cdot ( \hbox{the weight of $T$ on $W/H$}). $
  \end{enumerate}
\end{enumerate}
From these it is easy to see that the multidegree is a positive sum of
monomials in the weights of $T$ on $W$, and is homogeneous of degree
$\codim_W X$.
Also, for our action of $T$ on $\MNC$, 
$$ \deg_\MNC X = (\mdeg_\MNC X)|_{z_i \equiv 0, \A = 1}. $$
We will occasionally use a more general version of (3'b), in which $X$
and $Y$ are $T$-invariant varieties in $W$ whose intersection has the 
expected dimension. 
Then $\mdeg_W (X\cap Y) = (\mdeg_W X) (\mdeg_W Y)$. The most commonly
used case is when $Y$ is a $T$-invariant hypersurface, and $\mdeg Y$
is the weight of its defining equation. 

\newcommand\MNCdzero{\MNC_{\Delta=0}}

In this paper our interest is in the multidegrees of the $\{E_\pi\}$,
which all live in the subspace $\MNCdzero \leq \MNC$ of matrices
with zero diagonal. Hereafter we will drop the subscript on $\mdeg$, 
and assume it to be $\MNCdzero$,
with weights $\{\A + z_i - z_j\}, i,j \in\{1,2,\ldots,N\}, i\neq j$.

We can now state our main result (theorem \ref{thm:main}):
\begin{Theorem*}
  For all link patterns $\pi$, the
  multidegree of $E_\pi$ is the homogenized component $\Psi_\pi$:
  \begin{equation*}
    \mdeg E_\pi|_{\A=1}=\Psi_\pi
  \end{equation*}
\end{Theorem*}

The previous theorem is recovered by setting $z_i \equiv 0$.

We can give a geometric meaning to theorem 5 in \cite{DFZJ} (generalized
beyond $N$ even):

\begin{Corollary}\label{cor:totalmdeg}
  The multidegree of $E$ is the Pfaffian
  \begin{equation*}
    \mdeg E
    = {\rm Pf}\left({z_i-z_j\over \A-(z_i-z_j)^2}\right)_{1\leq i,j\leq N}
    \ \times \prod_{1\leq i<j\leq N} {\A-(z_i-z_j)^2\over z_i-z_j}
  \end{equation*}
  where ${\rm Pf}\ a={1\over 2^n n!}\sum_{\sigma\in S_N}
(-1)^\sigma \prod_{i=1}^n a_{\sigma(2i-1)\,\sigma(2i)}$.
%with $n=\lfloor N/2\rfloor$.

  Its ordinary degree is the determinant 
  $\det\Big[{2i+2j+1\choose 2i}\Big]_{0\leq i,j\leq n-1}=1,7,307,\ldots$
  for $N$ even, 
  $\det\Big[{2i+2j+3\choose 2i+1}\Big]_{0\leq i,j\leq n-1}=3,55,6153,\ldots$
  for $N$ odd.
\end{Corollary}

Since these multidegrees are very difficult to calculate directly, the
reader may wonder what insight has been gained about the $\{\Psi_\pi\}$.

\begin{Corollary}\label{cor:positive}
  Each coefficient $\Psi_\pi$ can be written as a sum, with positive
  coefficients, of products of distinct factors $\{1 + z_i - z_j\}$,
  where $i\neq j$ run over $1,\ldots,N$.
\end{Corollary}

\begin{proof}
  In general, if $T$ acts on $V$ with weights $\lambda_1,\ldots,\lambda_m$,
  the multidegree of a $T$-invariant subscheme $X \subseteq V$ 
  can be written as $p(\lambda_1,\ldots,\lambda_m)$,
  where $p$ is a polynomial with positive coefficients,
  each monomial squarefree.

  In the case at hand, $X=E_\pi$ and $V = \MNCdzero$,
  hence has weights $\{\A+z_i-z_j\}$ for $i\neq j$. 
\end{proof}

This positivity seems difficult to prove directly from equations
(\ref{eqn:pizeroa}) and (\ref{eqn:schub}), in much the same way that the
divided-difference definition of Schubert polynomials does not make it
easy to see that they have positive coefficients.

The most interesting case of the de Gier--Nienhuis conjecture was the
first observed: for $\pi(i)=N+1-i$, $N$ even, the coefficient
$d_\pi$ is the degree of the {\em commuting variety} of $n\times n$ matrices.
The recursion relations provide a formula, albeit rather involved, for
this quantity (theorem \ref{thm:comvar}):

\newcommand\MnC{M_n(\complexes)}
\begin{Theorem*}
  The degree of the commuting variety 
  $C=\{ (X,Y)\in \MnC : XY=YX\}$ is
  $$
  \deg C= \left[ \theta_1\ (\theta_2 \theta_1)\ \cdots\ (\theta_i
    \theta_{i-1}\cdots \theta_2 \theta_1) \ \cdots\ 
    (\theta_{n-1}\cdots\theta_2\theta_1) \prod_{i=1}^n
    (1+z_i)^{i-1}(1-z_i)^{n-i} \right] \Big|_{z_i\equiv 0}
  $$
  where $\theta_i=-2\der_i-\tau_i$. (For computational purposes, note that
  one can and should immediately specialize $z_i$ to $0$ after the
  last application of $\theta_{i-1}$, for each $i=1\ldots n$.)
\end{Theorem*}

In this way the connection is useful in the other direction.
In \cite{dGN} the authors used their Markov process to compute the degree of
the commuting variety (or at that point, a number conjecturally equal)
up through $8\times 8$ matrices.

Alternatively, one can use the formula above,
or rather, a slight simplication of it via a change of variable
proposed in \cite[section 6.2]{DFZJ}, as an efficient
algorithm for the computation of these numbers, which allows us to go
further. Here are the degrees through $11 \times 11$ matrices:
\begin{eqnarray*}
&\deg C=
1,\ 3,\ 31,\ 1145,\ 154881,\ 77899563,\ 147226330175,\ 1053765855157617,
 \\
&28736455088578690945,\ 3000127124463666294963283,\ 
1203831304687539089648950490463,\ldots
\end{eqnarray*}

\subsection{Acknowledgments}
We are thankful to Philippe Di Francesco, Edward Frenkel, Jan de Gier,
and Mark Haiman for useful conversations.

Throughout the paper, we use the notation 
$$ [P] = 
\begin{cases}
  1 &\mbox{if $P$ is true}\\
  0 &\mbox{if $P$ is false}\\
\end{cases}
$$
where $P$ is a property that may be true or false. For example, 
$\delta_{ab} = [a=b]$. We use $e^{ij}$ to indicate the matrix with
$1$ in entry $(i,j)$ and $0$ elsewhere, so
$$ (e^{ij})_{kl} = [i=k\hbox{ and }j=l].$$

The paper is organized as follows. Section 2 provides various definitions and properties
of the product $\cp$. Section 3 describes the Brauer loop scheme $E$ and its irreducible
components. Section 4 discusses their multidegrees and contains the main theorems of the
paper. The last four sections contain various additional results:
section 5 discusses the connection of the Brauer loop scheme with the previously
introduced upper-upper scheme and the application to the commuting variety, 
section 6 provides a geometric interpretation of some
recursion relations satisfied by the multidegrees, section 7 explains the connection
of the Brauer loop scheme to nilpotent orbits, and section 8 briefly mentions the
existence of a larger torus action.

\section{Models of $(\MNC,\cp)$}

Let $M_\leq$ denote the upper triangle of a matrix $M$, so 
$ (M_\leq)_{ij} = \begin{cases}
  M_{ij} &\mbox{if $i\leq j$}\\
  0 &\mbox{if $i>j$.}
\end{cases} 
$ \hfill\break
We will later use $M_>,M_<$ for the strict lower and upper triangles.

\subsection{The semidirect product model}\label{ssec:semidirect}

\newcommand\Upt{R_N(\complexes)} % upper triangular matrices
\newcommand\Utwopt{R_{2N}(\complexes)} % upper triangular matrices
\newcommand\UZ{R_\integers(\complexes)}
\newcommand\Uptaff{R_N(\complexes[t])} % affine upper triangular matrices

We can study the multiplication $(\MNC,\cp)$ in terms of ordinary
matrix multiplication, at the expense of making the cyclic invariance
less obvious. 

Let $\Upt$ denote the algebra of upper triangular matrices (on which the
usual product and the $\cp$ product coincide) and $\MNC/\Upt$ the evident 
quotient bimodule for $\Upt$.
Then the semidirect product $\Upt \times \MNC/\Upt$ carries the multiplication
$$ (R,L)\cdot (V,M) := (RV, RM+LV). $$
Our algebra $(\MNC,\cp)$ is isomorphic to this semidirect product,
via the map
\begin{eqnarray*}
  (\MNC,\cp) &\to& \Upt \times (\MNC/\Upt) \\
  M &\mapsto& (\hbox{$M_\leq$, $M + \Upt$}).
\end{eqnarray*}
An element $(R,L)$ of this semidirect product is invertible
(with inverse $(R^{-1}, -R^{-1} L R^{-1})$) if and only if $R$ is an
invertible upper triangular matrix, which of course is equivalent to
having all its diagonal elements be nonzero. There is no condition on $L$.

The group of units $(\MNC,\cp)^\times$ in this algebra
is therefore also a semidirect product 
$$ (\MNC,\cp)^\times \iso B \ltimes \MNC/\Upt, $$
where $B$ denotes the group of invertible upper triangular matrices,
and $\MNC/\Upt$ the vector space considered as an abelian group.
Hence $(\MNC,\cp)^\times$ is solvable, with the diagonal matrices
serving as a maximal torus, and the group $U$ as the unipotent radical.

\newcommand\Lie{{\rm Lie\,}}
\newcommand\Ad{{\rm Ad\,}}
\newcommand\ad{{\rm ad\,}}
In these $\Upt \times \MNC/\Upt$ coordinates, the scheme $E$ and the
action take the form
\begin{eqnarray*}
  E &\iso& \big\{ (R,L) : % U\in\Upt, L \in \MNC/\Upt, 
          R^2 = 0, \diag(R)=0, RL+LR \in \Upt \big\} \\
  \Ad (X,0)\cdot (R,L) &=& (XUX^{-1}, XLX^{-1}) \\
  \Ad (1,Y)\cdot (R,L) &=& (R, L + [Y,R])
\end{eqnarray*}
where we remember to always interpret the second entry as being in the
quotient space $\MNC/\Upt$.

\subsection{The periodic strip model}

Let $\UZ$ denote the space of upper triangular matrices $M$,
where the indices in $M_{ij}$ run over $\integers$. Despite the infinitude,
each sum
$$ (AB)_{ik} = \sum_j A_{ij} B_{jk}, \quad i\leq j\leq k $$
defining the product is finite.

Let $S \in \UZ$ denote the shift matrix with entries 
$S_{ij} = \delta_{i,j-1}$. 
To specify an element of the quotient ring $\UZ / \< S^N \> $,
one can use the matrix entries $L_{ij}$, $0 \leq j-i < N$, as
those with $0>j-i$ are zero by triangularity and with $j-i\geq N$ are
rendered ambiguous by the quotient. These ring elements can be
pictured as infinite diagonal strips of width $N$, on and above the
main diagonal.

\begin{Proposition}
  There is an injective ring homomomorphism
  $$ \Phi : (\MNC,\cp) \to \UZ / \< S^N \> $$
  given by $\Phi(M)_{ij} = M_{i\bmod N,j\bmod N}$ for $0\leq j-i<N$.
  The image is the space of periodic strips, i.e.
  $\Phi(M)_{ij} = \Phi(M)_{i+N,j+N}\ \forall 0\leq j-i<N$.
\end{Proposition}

\begin{proof}
  The only claim worth commenting on is the ring homomorphism.
  Let $i,k$ satisfy $0\leq k-i < N$. Since $\Phi(M)$ is periodic, 
  we will assume $1\leq i\leq N$ as well. Then there are two cases,
  depending on $k\leq N$ or $k>N$. If $k\leq N$,
  \begin{align*}
  \Phi(P\cp Q)_{ik} &= (P\cp Q)_{ik} 
  = \sum_{j, i\leq j\leq k} P_{ij} Q_{jk} 
  = \sum_{j, i\leq j\leq k} \Phi(P)_{ij} \Phi(Q)_{jk} \\
  &= (\Phi(P) \Phi(Q))_{ik}. 
  \end{align*}
  Whereas if $k>N$,
  \begin{eqnarray*}
  \Phi(P\cp Q)_{ik} 
  &=& (P\cp Q)_{i,k-N} 
  = \sum_{j, i\leq j\leq N} P_{ij} Q_{j,k-N} 
          + \sum_{j, 1\leq j\leq k-N} P_{ij} Q_{j,k-N} \\
  &=& \sum_{j, i\leq j\leq N} \Phi(P)_{ij} \Phi(Q)_{jk} 
          + \sum_{j, 1\leq j\leq k-N} \Phi(P)_{i,j+N} \Phi(Q)_{j+N,k} \\
  &=& \sum_{j, i\leq j\leq k} \Phi(P)_{ij} \Phi(Q)_{jk} 
  = (\Phi(P) \Phi(Q))_{ik}. 
  \end{eqnarray*}
\end{proof}

This model has the benefit of making the cyclic invariance obvious, and
is the easiest to calculate with visually. To connect it with the 
semidirect product model, a pair $(R,L)$ corresponds to the strip
$$ \begin{array}{ccccc}
  \ddots&L&&& \\
  &R&L&& \\
  &&R&L& \\
  &&&R&L \\
  &&&&\ddots
  \end{array} 
$$

Such representations as periodic infinite matrices have been considered
in the context of loop algebras \cite{KR}. This leads us naturally to
the next model:

\subsection{The affine $GL_n$ model}\label{ssec:algdegen}

Consider the ring of matrices $\MNCt$ over the polynomial ring $\complexes[t]$.
Let $\Uptaff$ denote the $\complexes[t]$-subalgebra generated by
$\{e^{i,i+1}\}$ for $1\leq i<N$, and $t\, e^{N,1}$. (These are the simple
root spaces of the affine Lie algebra $\widehat{\lie{gl}_N}$, leading to the
name of this model.) Then the
following is straightforward from the semidirect product model:

\begin{Proposition}
  The map $(\MNC,\cp) \to \Uptaff/(t\Uptaff)$ taking $M$ to 
  $M_\leq + t M_>$ is an isomorphism.
\end{Proposition}

We can regard $\Uptaff$ as a family of algebra structures on $\MNC$ 
indexed by $t \in \complexes$,
where the fiber $t=1$ is ordinary multiplication and $t=0$ is $\cp$.
There is an associated flat family whose $t$-fiber is the space of
matrices that square to zero under the $t$-multiplication. 
We investigate this family in section \ref{sec:flatlimit}, where we show
that the flat limit as $t\to 0$ is supported on the top-dimensional components 
of $E$, and contains each component with multiplicity $2^{\lceil N/2\rceil}$. 

There is another way, manifestly cyclically invariant, to degenerate the
algebra $(\MNC,\times)$ to the algebra $(\MNC,\cp)$. Let $s\cdot M$ be
defined by
$$ (s\cdot M)_{ij} = s^{(j-i)\bmod N} M_{ij}, \qquad (j-i) \bmod N \in [0,N) $$
and define 
$$ M\times_s N := s^{-1} \cdot ((s\cdot M) (s\cdot N)). $$
So $M \times_1 N = MN$, and for $s\neq 0$ this multiplication is conjugate
to the ordinary one. (If we left out the ``$\bmod\ N$'' part, it would be
{\em equal} to the ordinary one.) Then it is easy to check that 
$$ \lim_{s\to 0}\ M\times_s N = M\cp N. $$

\section{Components of $E$ and link patterns}

\subsection{Decomposition of $E$ in terms of involutions}

Recall that we use $\Mpi$ to denote the permutation matrix of a
permutation $\pi$, with the diagonal zeroed out.  We care especially
about involutions, because of Melnikov's theorem:

\begin{Theorem}\cite{M}\label{thm:melnikov}
  Let $B := \Upt^\times$ denote the group of $N\times N$ invertible
  upper triangular matrices.  The action by conjugation of $B$ on the
  set $\{X \in \Upt : X^2 = 0\}$ has finitely many orbits, and each
  contains a unique partial permutation matrix.

%  Moreover, if $X = Ut$ where $t$ is diagonal, then $X$ can be conjugated
%  to $U$ using a diagonal matrix.

  A partial permutation matrix is an element of this space if and only if
  it is $\Mpisl$ for some involution $\pi \in S_n$.
  Hence the orbits are naturally indexed by involutions.
\end{Theorem}

For example, the identity matrix is an involution whose strict upper
triangle vanishes, and the corresponding orbit consists only of the
zero matrix. Using the semidirect product model, we easily obtain

\begin{Corollary}\label{cor:Rinv}
  If $M\in E$, then there exists a $\cp$-invertible $P$ and
  an involution $\pi$ such that 
  $$ (P\cp M \cp P^{\cp-1})_\leq = \Mpisl.$$
\end{Corollary}

\begin{Theorem}\label{thm:diagpaired}
  Let $M\in E$. Then the ordinary square $M^2$ (not $M\cp M$) has 
  diagonal entries which come in pairs, or put another way, there exists 
  a link pattern $\pi$ such that $(M^2)_{ii} = (M^2)_{\pi(i)\pi(i)}$
  for each $i=1,\ldots,N$.
  
  Conversely, every link pattern $\pi$ is necessary: 
  there exists $M\in E$ such that $(M^2)_{ii} = (M^2)_{jj}$ if and
  only if $j\in \{i,\pi(i)\}$.
\end{Theorem}

\begin{proof}
  We first check that these diagonal elements are invariant under conjugation.
  If $M'=P\cp M\cp P^{\cp-1}$,
  \begin{equation*}
    (M'^2)_{ii}=\sum_{l=1}^N M'_{il} M'_{li}
    = \sum_{\scriptstyle j,k,l,p,q\atop
      {\scriptstyle \cyc{i\le j\le k\le l}\atop
        \scriptstyle \cyc{l\le p\le q\le i}}}
    P_{ij} M_{jk} P^{\cp-1}_{kl}
    P_{lp} M_{pq} P^{\cp-1}_{qi}
  \end{equation*}
  
  This can be visualized with $i,j,k,l,p,q,i$ winding only once 
  counterclockwise round
  a circle. If $p=k\ne l$, then $p=q=i=j=k$, and these terms contain a
  factor $M_{ii}^2=0$. The remaining terms have $\cyc{k\le l\le p}$, 
  so that one can perform the summation over $l$:
  \begin{equation*}
 (M'^2)_{ii}=\sum_{\scriptstyle j,k,q\atop 
 {\scriptstyle \cyc{i\le j< k}
 \atop\scriptstyle \cyc{k< q\le i}}}
  P_{ij} M_{jk} M_{kq} P^{\cp-1}_{qi}
  \end{equation*}
  Let us consider the summation at fixed $j$ and $q$. If $i\ne j$ or
  $q\ne i$, one finds $\cyc{j\le k\le q}$, and the sum over $k$ is equal
  to $(M\cp M)_{jq}$, which is zero for $M\in E$.  There remain only
  the contributions at $i=j$ and $q=i$, which reduce to
  $(M'^2)_{ii}=(M^2)_{ii}$.

  Now we use corollary \ref{cor:Rinv} to reduce to the case that 
  $M_\leq = \Mpisl$ for some involution $\pi$ (not necessarily a 
  link pattern). Then one easily computes
  $$ (M^2)_{ii} = 
  \begin{cases}
  M_{i,\pi(i)} & \mbox{if $i> \pi(i)$} \\
  M_{\pi(i),i} & \mbox{if $i< \pi(i)$} \\
  0 & \mbox{if $i = \pi(i)$} 
  \end{cases}
  $$
  and hence $(M^2)_{ii} = (M^2)_{\pi(i)\pi(i)}$.
  
  To see that every link pattern $\pi$ arises, let $t$ be a generic
  diagonal matrix, and $M = \Mpi t$.  Then
  $$ (M^2)_{ii} = ((\Mpi t)^2)_{ii} = 
  \begin{cases}
    t_i t_{\pi(i)} & \mbox{ if $i\neq \pi(i)$} \\ \\
    0 & \mbox{ if $i = \pi(i)$.} 
  \end{cases} 
  $$
  By the genericity, $t_i t_{\pi(i)} \neq t_j t_{\pi(j)}$ unless
  $i=j$ or $i=\pi(j)$, and $t_i t_{\pi(i)} \neq 0$. Since $\pi$ is a
  link pattern, there is at most one $0$. So the only repetitions are
  the expected ones.
\end{proof}

\begin{Theorem}\label{thm:topcomps}
  The scheme $E$ is $\lfloor N^2/2\rfloor$-dimensional, and the top
  components correspond to link patterns. Moreover, for each link
  pattern $\pi$ the scheme
  $$ E_\pi := \overline{
    \big\{ M\in E 
    : (M^2)_{ii} = (M^2)_{jj}\hbox{ if and only if } j\in\{i,\pi(i)\} \big\}}$$
  is irreducible.
\end{Theorem}

\begin{proof}
  We will give a finite decomposition of $E$ into irreducible pieces 
  $\{F_\pi\}$ corresponding to involutions, with 
  $\dim F_\pi = \frac{1}{2}(N^2$ minus the number of fixed points of $\pi$). 
  The closures of the $F_\pi$ of largest dimension are definitely components,
  and there may be other, smaller components.\footnote{Brian Rothbach
    has shown there are not.}

  Consider the map $\rho: E \to \Upt$ given by $(R,L)\mapsto R$,
  in the semidirect product model. If we let $B$ act on $E$ by
  $\cp$-conjugation (where we identify $B \iso \{ (R,0) : R$ invertible$\}$),
  and on $\Upt$ by ordinary conjugation,
  then this map $\rho$ is $B$-equivariant. By theorem \ref{thm:melnikov},
  the image is a finite union of $B$-orbits, with the set $\{\Mpisl : \pi$ 
  an involution$\}$ serving as orbit representatives.

  For $\pi$ an involution, let $F_\pi := \rho^{-1}(B \cdot \Mpisl)$,
  so $E$ is the finite disjoint union of the locally closed pieces $\{F_\pi\}$.
  Then restricted to $F_\pi$, the map $\rho$ is a fiber bundle
  (since the image is a $B$-orbit), and it is enough to understand one fiber. 
  In particular, 
  $$ \dim F_\pi = \dim (B\cdot \Mpisl) 
  + \dim\ \{L : L\Mpisl + \Mpisl L \in \Upt \}. $$
  The dimension of the $B$-orbit was computed in \cite[section 3.1]{M}
  (where it is called $m+s$), but we will not make direct use of the
  slightly intricate formula given there.

  Let $\ad X\cdot Y := XY-YX$. 
  Consider the map $\ad \Mpisl$ on $\Upt$
  \begin{equation}    \label{eq:adpi}
     e^{ij} \mapsto e^{i\pi(j)}\ [j<\pi(j)] - e^{\pi(i)j}\ [\pi(i)<i], 
  \qquad i<j 
  \end{equation}
  whose image is the tangent space to $B\cdot \Mpisl$ at $\Mpisl$. 
  It contains the subspace $(\ad \Mpisl) \cdot \Upt_+$ where $\Upt_+$
  is the {\em strictly} upper triangular matrices. This subspace
  $(\ad \Mpisl)\cdot \Upt_+$ has codimension $m$ in $(\ad \Mpisl)\cdot \Upt$,
  where $m$ is the number of $2$-cycles in $\pi$, as easily seen by
  applying $\ad \Mpisl$ to the diagonal matrices.

  Now consider the equations $\{L : \Mpisl L + L\Mpisl \in \Upt \}$.
  For each $i<j$, the lower triangle entry $(j,i)$ must vanish:
  $$ (\Mpisl L + L\Mpisl)_{ji} 
  = L_{\pi(j)i}\ [j<\pi(j)] + L_{j\pi(i)}\ [\pi(i)<i] = 0, \qquad i<j $$
  Let $M_{ab} = L_{ab}$, times $-1$ if $\pi(b)<b$. Then these 
  restrictions on $L$ are equivalent to 
  \begin{equation}    \label{eq:Mconditions}
    M_{\pi(j)i}\ [j<\pi(j)] - M_{j\pi(i)}\ [\pi(i)<i] = 0, \qquad i<j. 
  \end{equation}
  (The signs only matter when both terms appear, and in this case it
  is easy to check that only the second one is negated.) 

  Since (\ref{eq:adpi}) and (\ref{eq:Mconditions}) have the same form,
  the space of matrices $\{M \in \MNC/\Upt\}$ satisfying these conditions
  (\ref{eq:Mconditions})
  is exactly the perpendicular to the space $(\ad \Mpisl)\cdot \Upt_+$
  spanned by the image of (\ref{eq:adpi}), where ``perpendicular'' is
  defined with respect to
  the perfect pairing $\langle R,L\rangle := \Tr (RL)$ between
  $\MNC/\Upt$ and $\Upt_+$. Hence
  $$ \dim\ (\ad \Mpisl)\cdot \Upt_+ 
  + \dim\ \{L : \Mpisl L + L\Mpisl \in \Upt \} = \dim \Upt_+ 
  = \onehalf (N^2-N). $$
  With $m$ more from $(\ad \Mpisl)\cdot \Upt / (\ad \Mpisl)\cdot \Upt_+$,
  the dimension of $F_\pi$ is $\onehalf (N^2-N)+m$.
  This is only maximized when $m=n$, i.e. $\pi$ is
  a link pattern.

  Finally, since $F_\pi$ is a fiber bundle over the $B$-orbit 
  $B\cdot \Mpisl$ with fiber a vector space \break
  $\{L : \Mpisl L + L\Mpisl \in \Upt \}$, it is
  irreducible. So each $\overline{F_\pi}$, for $\pi$ a link pattern,
  contributes only one component of top dimension to $E$.

  By the computation at the end of theorem \ref{thm:diagpaired}, 
  the set $\big\{ M\in E : (M^2)_{ii} = (M^2)_{jj}$ if and only if
  $j\in\{i,\pi(i)\} \big\}$ is contained in $F_\pi$. Since they have
  the same dimension, this subset too is irreducible, 
  as is its closure $E_\pi$.
\end{proof}

A similar technique was used in \cite[lemma 1]{Kn} to determine the
components of the upper-upper scheme. 
In that case the 
dimension bound lets one prove that the upper-upper scheme is a
complete intersection, hence has no lower-dimensional components. 
Brian Rothbach has shown this equidimensionality also holds for the Brauer
loop scheme $E$ (which is not a complete intersection).

\subsection{Properties of the $\{E_\pi\}$ components}

In this section we show that the components $\{E_\pi\}$ are
generically reduced, we parametrize them, and find some (and
conjecturally, all) of their defining equations.

\begin{Theorem}\label{thm:reduced}
  Each $E_\pi$ is reduced at $\Mpi t$ for $t$ generic diagonal.
  Hence $E_\pi$ is generically reduced.
\end{Theorem}

\begin{proof}
%  It is enough to prove reducedness at one point. 
  We do this by showing
  that the Zariski tangent space has the right dimension.
  The Zariski tangent space is the common kernel of the derivatives at
  $\Mpi t$ of the defining equations for $E$.
  
  The linear equations are handled by just working inside the
  $(N^2-N)$-dimensional space $\MNCdzero$. The derivative of 
  $M\cp M = 0$ is $P \mapsto P\cp M + M \cp P$. In the case at hand, 
  \begin{eqnarray*}\label{eqn:zariski}
    (P\cp (\Mpi t)+(\Mpi t)\cp P)_{ik} 
    &=& \sum_{j:\ \cyc{i,j,k}} (P_{ij} (\Mpi t)_{jk} + (\Mpi t)_{ij} P_{jk}) \\
    &=& P_{i \pi(k)} t_k [\cyc{i\le\pi(k)\le k}]
    + t_{\pi(i)} P_{\pi(i) k} [\cyc{i\le\pi(i)\le k}].
  \end{eqnarray*}

  We require these to be zero for all $i$ and $k$.
  Let us organize the equations as follows. If $i=k$ or $\pi(k)$ the
  equation is trivial.  So we can assume that $i$ and $k$ belong to
  distinct orbits.  Diagramatically, there are three ways for
  the orbits $\{i,\pi(i)\}$, $\{k,\pi(k)\}$ to relate:
  
  \begin{enumerate}
  \item The chords $\{ i,\pi(i)\}$ and $\{ k,\pi(k)\}$ cross each other.
    In this case we can choose the labelling so that
    $\cyc{i< k<\pi(i)<\pi(k)}$:   
    \begin{minipage}{2cm}\epsfig{file=cross.eps, width=2cm}\end{minipage},
    and by inspection we find the
    following four equations:
    \begin{eqnarray*}
      t_i P_{ik} + t_k P_{\pi(i)\pi(k)} &=& 0\\
      t_k P_{k\pi(i)} + t_{\pi(i)} P_{\pi(k)i} &=&0 \\
      t_{\pi(i)} P_{\pi(i)\pi(k)} + t_{\pi(k)} P_{ik} &=&0 \\
      t_{\pi(k)} P_{\pi(k)i} + t_i P_{k\pi(i)} &=&0
    \end{eqnarray*}
    (all these equations are obtained from each other by rotation of
    $90^\circ$, which is the symmetry of the diagram). Generically, 
    $t_i t_{\pi(i)}\ne t_k t_{\pi(k)}$ and we can in fact simplify this
    system to
    \begin{eqnarray*}
      P_{ik}=
      P_{k\pi(i)}=
      P_{\pi(i)\pi(k)}=
      P_{\pi(k)i}=0
    \end{eqnarray*}
    which shows that there are exactly four independent equations.
  \item The chords $\{i,\pi(i)\}$ and $\{k,\pi(k)\}$ do not cross each
    other, in which case we can choose $\cyc{i<\pi(i)< k<\pi(k)}$:
    \begin{minipage}{2cm}\epsfig{file=noncross.eps, width=2cm}\end{minipage}.
    We find again four equations, though of a different form:
    \begin{eqnarray*}
      t_i P_{ik} + t_k P_{\pi(i)\pi(k)} &=& 0 \\
      t_k P_{ki} + t_i P_{\pi(k)\pi(i)} &=& 0\\
      t_{\pi(i)} P_{\pi(i)k} &=& 0 \\
      t_{\pi(k)} P_{\pi(k)i} &=& 0
    \end{eqnarray*}
    (note that they form groups of two, related by a rotation of
    $180^\circ$ or equivalently exchange of $i$ and $k$).  They are
    generically (for non-zero $t$'s) non-trivial and independent from
    each other.
  \item If one of the indices is a fixed point, one can assume that
    $\cyc{i<\pi(i)< k=\pi(k)}$, in which case one finds two equations:
    \begin{eqnarray*}
      P_{\pi(i)k}=
      P_{ki}=0
    \end{eqnarray*}
  \end{enumerate}
  (Not both $i$ and $k$ can be fixed, since $i\neq k$ and $\pi$ is
  a link pattern.)
  
  The conclusion is that each pair of chords contributes exactly $4$
  equations, and a chord plus a fixed point contributes $2$ equations;
  thus, recalling that $N=2n+r$ with $r=0,1$ the number of fixed points,
  a total of $4\times n(n-1)/2+2\times n r=2n(n+r-1)$
  equations.  Therefore the kernel is of dimension
  $N(N-1)-2n(n+r-1)=2n(n+r)+r(r-1)$. Setting $r=0$, $1$ we find the
  desired dimension $2n(n+r) = \lfloor N^2/2\rfloor$.
\end{proof}

At this point we have three equivalent definitions of $E_\pi$:
\begin{itemize}
\item the closure of $\{M \in E : (M^2)_{ii} = (M^2)_{jj} \Longleftrightarrow
  j \in \{i,\pi(i)\} \ \}$
\item the closure of $\{M \in E : M_<$ is $B$-conjugate to $\Mpisl\}$
\item the unique component of dimension $\lfloor N^2/2\rfloor$ in
  $\{M \in E : (M^2)_{ii} = (M^2)_{\pi(i)\pi(i)} \}$.
\end{itemize}

This third definition is a first step in defining $E_\pi$ by equations.
To do better, we use yet another characterization of $E_\pi$.

\begin{Proposition}\label{prop:openset}
  Let $\pi$ be a link pattern, and $\Mpi$ its permutation matrix with
  the diagonal zeroed out.
  The irreducible set $U\cdot\{\Mpi t, t\in T\}$ is dense in $E_\pi$.
\end{Proposition}

\begin{proof}
  Since $U$ and $T$ are irreducible, so is $U\cdot \{\Mpi t, t\in T\}$.
  Following the calculation at the end of theorem \ref{thm:diagpaired},
  we see $U\cdot \{\Mpi t, t\in T\} \subseteq E_\pi$.
  
  There are two steps. The first is to compute the dimension of 
  a generic $U$-orbit $U\cdot (\Mpi t)$.
  The second is to show that each $U$-orbit intersects
  the set of representatives $\{\Mpi t\}$ in only one point,
  hence the dimension of $U\cdot\{\Mpi t, t\in T\}$ 
  is the dimension of $\{\Mpi t, t\in T\}$ plus the
  dimension of a generic $U$-orbit.
% (In fact it would be enough to
%  show it intersects in finitely many points, but our argument leads to $1$.)

  We now compute the infinitesimal stabilizer of $U$ on $\Mpi t$, 
  where $t$ is generic.
  Let $P$ be an element of the Lie algebra of $U$, which is $\MNCdzero$.
  The equation $\Mpi t\cp P=P\cp \Mpi t$ reads
  \begin{equation*} \label{eqn:proofdim}
    t_{\pi(i)} P_{\pi(i)k} [\ \cyc{i\le\pi(i)\le k}] 
    = P_{i\pi(k)}t_k [\ \cyc{i\le\pi(k)\le k}].
  \end{equation*}
  
  Note that these equations are exactly of the same form as those
  %equation \ref{eqn:zariski} 
  in the proof of theorem~\ref{thm:reduced}, up to a
  sign (much as went into equation (\ref{eq:Mconditions})), and we
  shall not repeat the arguments that lead to the conclusion that 
  $P$ satisfies $2n(n+r-1)$ equations %with $N=2n+r$,  $r=0$ or $1$, 
  and therefore this is also the dimension of $U\cdot (\Mpi t)$.
  
  Next, assume that $P\cp \Mpi t=\Mpi t'\cp P$. For each $i=1,\ldots,N$,
  the equation concerning entry $(i,\pi(i))$ reads $t_i=t'_i$.
  So each $U$-orbit contains a unique element of the form $\Mpi t$.

  Finally, noting that $\dim \{\Mpi t,t\in T\}=2n$,
  we compute $\dim U\cdot\{\Mpi t, t\in T\}=2n(n+r)=\lfloor N^2/2\rfloor$.
  Since $U\cdot\{\Mpi t, t\in T\} \subseteq E_\pi$ and has the same 
  dimension, it is dense in $E_\pi$.
\end{proof}

Any equations satisfied by this dense open set are satisfied by all of $E_\pi$.
We pay special attention to the linear equations, mostly in order
to connect to proposition~1 of \cite{DFZJ}. 

\begin{Proposition}\label{prop:linearcondition}
  Assume the link pattern $\pi$ has no chord connecting a pair of
  points between labels $i$ and $l$ (i.e.\ there are no $j$ s.t.\ 
  $\cyc{i\le j\le l}, \cyc{i\le\pi(j)\le l}$). Then $M\in E_\pi$ implies
  $M_{il}=0$. 
  
  More generally, the periodic strip associated to $M$ vanishes %everywhere
  southwest of the $(i,l)$ entry.
\end{Proposition}

\begin{proof} 
  By the density, it is enough to check for $M = P\cp (\Mpi t)\cp P^{\cp-1}$ 
  for some diagonal $t$.  Write
  \begin{equation*}
    M_{il}=\sum_{\cyc{i, j, k, l}} P_{ij}\,  \Mpi_{jk} t_k {P^{\cp-1}}_{kl}
  \end{equation*}
  and notice that $k=\pi(j)$, $\cyc{i, j, k, l}$ contradicts the
  hypothesis on $\pi$.  Therefore the sum is zero.
  
  For the second conclusion, note that the hypothesis for the pair
  $(i,l)$, plus $\cyc{i\le j\le k\le l}$, implies the hypothesis for
  the pair $(j,k)$.
\end{proof}

Let $r_{ij}(M)$ denote the rank of the triangular matrix southwest of
the $(i,j)$ entry in the periodic strip model of $M$. In this language,
the previous proposition asserted that $r_{il}(M) = 0$ for certain $(i,l)$.

\begin{Theorem}\label{thm:compeqns}
  The variety $E_\pi$ satisfies the following equations:
  \begin{enumerate}
  \item those defining $E:\ M\cp M = 0$
  \item those defining $E_\pi :\ (M^2)_{ii} = (M^2)_{\pi(i)\pi(i)}$
  \item for any $M \in E_\pi$, and matrix entry $(i,j)$,
    we have $r_{ij}(M) \leq r_{ij}(\Mpi)$. In polynomial terms, this 
    asserts the vanishing of all the minors of size $r_{ij}(\Mpi)+1$
    in the submatrix southwest of entry $(i,j)$ in the strip model.
  \end{enumerate}
\end{Theorem}

\begin{proof}
  The first two are automatic. 
  For the third, note that the action of $U$ in the periodic strip model
  has a well-defined restriction to each southwest triangle, since
  $U$ acts by north- and east-moving row and column operations.
\end{proof}

\begin{figure}[htbp]
  \centering
  \epsfig{file=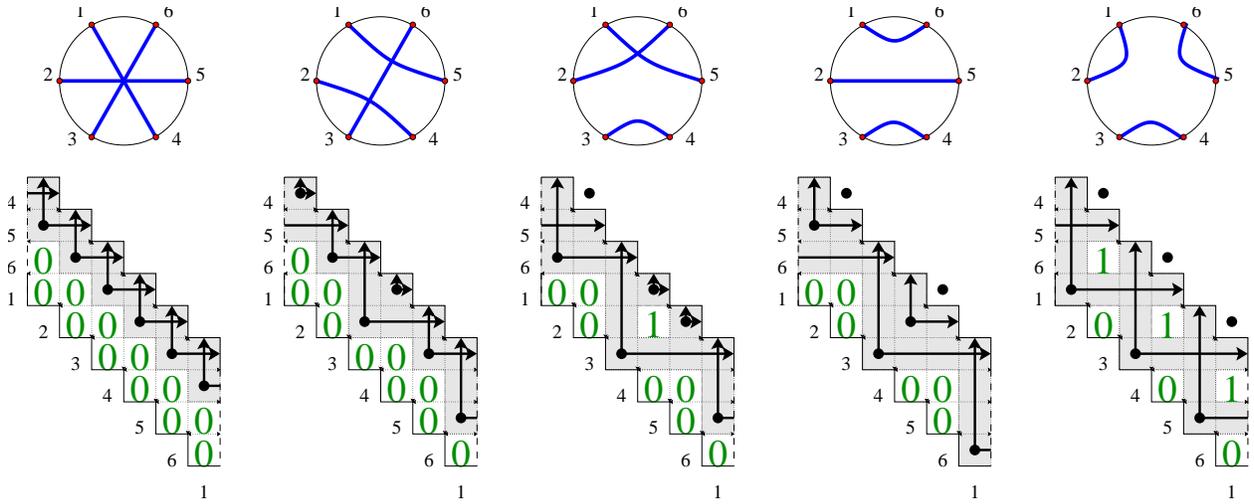,width=6.5in}
  \caption{The diagrams for the link patterns of size $6$, up to rotation.
    (We have left out the top diagonal $\{(i,i-1)\}$, whose irrelevance
    will be shown in lemma \ref{lem:freeentries}.)  The rank function
    $r_{ij}$ is written in the boxes of each diagram; for example 
    $r_{ij}=0$ implies that that matrix entry is actually $0$.}
% For example, 
%    each of the rank $1$ conditions requires that the $2\times 2$
%    matrix to the southwest must be singular.}
  \label{fig:ess6}
\end{figure}

The third group of conditions appear in a similar context in \cite{Fu}, 
defining {\em matrix Schubert varieties}. These conditions are highly
interdependent, and Fulton defined the {\em essential set} $\{(i,j)\}$
whose rank conditions imply all the others. In the context at hand, 
the analogous set is defined as follows. Draw $\Mpi$ in the strip model,
and cross out all the boxes $(i,j)$ directly north or directly east of 
each $1$ entry in $\Mpi$. The remaining set of boxes in the strip is the 
{\dfn diagram} of the link pattern $\pi$, and the northeast corners
of each component of the diagram are the {\dfn essential set} of the
diagram. Then it is easy to check that the rank conditions 
$r_{ij}(M) \leq r_{ij}(\Mpi)$ for $(i,j)$ not in the essential set are
implied by those from the essential set. See figure \ref{fig:ess6} for
the possible diagrams in $N=6$.

Much the same sort of parametrization, and consequence for the
defining equations, was used in \cite{Kn} for the upper-upper scheme.

We conjecture that the equations in theorem \ref{thm:compeqns}
are {\em all} the equations defining the components.
Because of the connection explained in 
section \ref{sec:perm} between $E$ and the upper-upper scheme, 
this would imply a similar conjecture in \cite{Kn}, which in turn would
imply the well-known conjecture that the commuting scheme is reduced.

\section{Multidegrees and the Brauer loop model}\label{sec:brauer}

Our goal in this section is the main result of the paper, the equality 
$(\mdeg E_\pi)|_{\A=1} = \Psi_\pi$ relating $E$ to the Brauer loop polynomials.
The base case is easy: 

\begin{Proposition}\label{prop:pizero}
  %Set $n=\lfloor N/2 \rfloor$. 
  Define $\pi_0(i)=i+n \bmod 2n$ for
  $i\leq 2n$, and $\pi_0(N)=N$ if N is odd. 
  Then
  \begin{equation} \label{eqn:pizero}
    \mdeg E_{\pi_0} 
    =\prod_{i=1\ldots N \atop j:\ \cyc{i<j<i+n}} \!\!\!\!\!\! (\A+z_i-z_j)
  \end{equation}
  for $N$ even and
  \begin{equation*}
    \mdeg E_{\pi_0} 
    =\prod_{i=1\ldots N \atop j:\ \cyc{i<j<i+n}} \!\!\!\!\!\! (\A+z_i-z_j)\ \ 
    \left( \prod_{i=n+1}^N (\A+z_i-z_{i+n}) \right)
  \end{equation*}
  for $N$ odd.
\end{Proposition}

\begin{proof}
  Proposition \ref{prop:linearcondition} implies that $M_{il}=0$,
  $l=i,\ldots,i+n-1\mod N$ for $N$ even and similarly for $N$ odd.
  These linear equations imply $M\cp M=0$, and are of the right
  codimension ($2n(n-1)$ for $N$ even, $2n^2$ for $N$ odd) to define a
  component of $E$, so they exactly define $E_{\pi_0}$.  The result
  follows from property (3'b) and property (1) for multidegrees.
\end{proof}

\subsection{Geometry of divided difference operators}
The recursion relation (\ref{eqn:schub})
involves a divided difference operator $\der_i$,
so we take a moment to recall the geometry behind these,
making no claims to originality.

%We recall now one geometric derivation of divided difference operators. 
Let $L_i$ (for Levi subgroup) denote the group of invertible matrices with
off-diagonal elements only in entries $(i,i+1), (i+1,i)$. 
Let $B_i$ denote the subgroup in which the $(i+1,i)$ entry vanishes
(so $B_i\leq B$, unless $i=N$). 

Let $X$ carry a left $B_i$-action, let $V$ be a $L_i$-representation,
and let $f: X \to V$ be $B_i$-equivariant.
%(Soon $f$ will be the inclusion of $E^\iip_\pi$ into zero-diagonal matrices.)
(In our case $f$ will be an inclusion.)
Define the {\em map} $-\der_i f$ by
\begin{eqnarray*}
  -\der_i f : 
  L_i \times_{B_i} X &\to& V \\
  (P,M) &\mapsto& P \cdot f(M)
\end{eqnarray*}
where $L_i \times_{B_i} X$ is the quotient of $L_i \times X$ by the
diagonal action of $B_i$ on the right of $L_i$ and the left of $X$.
(Note that the formula stated descends to this quotient.)
One way to view $L_i \times_{B_i} X$ is as the associated $X$-bundle
to the principal $B_i$-bundle over $L_i/B_i \iso \Pone$. 

We now explain why this construction, based on the one of Bott and
Samelson \cite{BS}, Demazure, and Hansen, is given such a suggestive name:

\newcommand\Ker{{\rm Ker\,}}
\newcommand\Img{{\rm Im\,}}
\newcommand\Imgc[1]{\overline{\rm Im\, #1}}

\begin{Lemma}\label{lem:divdiff}
  Let $X$ be a variety in $V$ invariant under $B_i$ and rescaling, 
  with multidegree $\mdeg_V X$. Then
  $$ (-\der_i f)_*(1) = -\der_i \ \mdeg_V X $$
  where $(-\der_i f)_*$ denotes the pushforward map in equivariant cohomology.

  If $-\der_i f$ is generically $1$:$1$, then
  $$ (-\der_i f)_*(1) = \mdeg_V\ \Imgc{-\der_i f},$$
  the multidegree of the closure of the image.
\end{Lemma}

\begin{proof}
  This proof requires more equivariant cohomology than just the multidegree
  technology used elsewhere in the paper. There are many references,
  such as \cite{Br}.

  The space $L_i \times_{B_i} X$ equivariantly retracts to $L_i/B_i$ 
  (since it is an $X$-bundle, and $X$ can be rescaled toward $0$), 
  hence is equivariantly formal. So the map 
  $$ H^*_T(L_i \times_{B_i} X) \to H^*_T((L_i \times_{B_i} X)^T) $$
  is an inclusion. This factors through
  $$ H^*_T(L_i \times_{B_i} X) 
  \to H^*_T((B_i\times_{B_i} X)\ \cup\ (s_i B_i\times_{B_i} X)) $$
  since every $T$-fixed point in $L_i \times_{B_i} X$ lies over one of
  the two $T$-fixed points in $L_i/B_i$. Here $s_i$ denotes the permutation
  matrix of the transposition $(i \leftrightarrow i+1)$.
  
  Let $\alpha = z_i - z_{i+1}$. This is the weight on the tangent
  space $T_{s_i B_i/B_i}(L_i/B_i)$, and the negative of the weight on
  $T_{B_i/B_i}(L_i/B_i)$.  Then we claim the following equality of
  elements of $H^*_T(L_i \times_{B_i} X) \tensor_{H^*_T}   H^*_T[\alpha^{-1}]$,
  a localization of equivariant cohomology:
  $$ 1 
  = \frac{[B_i \times_{B_i} X]}{-\alpha} 
  + \frac{[s_i B_i\times_{B_i} X]}{\alpha} $$
  This is proved by restricting to the two $T$-invariant fibers
  $B_i \times_{B_i} X$ and $s_i B_i\times_{B_i} X$. Being the preimages of the
  points $B_i/B_i$ and $s_i B_i/B_i$, their normal bundles are trivial, with
  equivariant Euler classes $-\alpha$ and $\alpha$.

  When we apply $(-\der_i f)_*$ to both sides of this equation, we get
  $$ (-\der_i f)_*(1) = -\der_i (f_*(1)) = -\der_i \ \mdeg_V X, $$
  as desired.

  The second statement is standard.
\end{proof}

There are two differences between the usual Bott--Samelson construction and
the one used here. One is that Bott--Samelson maps usually take values in
a flag manifold, rather than a vector space, but this is immaterial to
the argument. The important one is that the input map $f$ is traditionally
equivariant under a {\em lower} triangular group $B_-$, so the tangent space
to the basepoint $B_- \in G/B_-$ has weights that are {\em positive} roots.
In our case we have invariance under the {\em upper} triangular $B_i$, 
leading us to the negative of the usual $\der_i$.

\newcommand\iip{{[i,i+1]}}
\subsection{The application to $\{E^\iip_\pi\}$}

\begin{Lemma}\label{lem:freeentries}
  Let $M \in \MNCdzero$, $\lambda \in \complexes$, 
  and $M' = M + \lambda e^{i+1,i}$. Then $M\cp M = M' \cp M'$, and in
  particular, $M\in E$ if and only if $M' \in E$.
\end{Lemma}

\begin{proof}
  If $e^{jk} \cp e^{i+1,i} \neq 0$, then $k=i+1$ and $\cyc{j\leq i+1\leq i}$,
  which forces $j=i+1=k$. Similarly 
  $e^{i+1,i} \cp e^{jk} \neq 0$ implies $j=k$. Since $M$ has zero diagonal 
  there are no such terms to contribute to the square.
\end{proof}

This lemma \ref{lem:freeentries} says that the entries on the top
diagonal in the strip model are unconstrained. As such, we don't lose
any information by setting them to $0$, and we learn something about
$\cp$-conjugation:

\begin{Lemma}\label{lem:conjactions}
  Let $\pi$ be a link pattern. Define 
  $$ E^\iip_\pi :=\{ M \in E_\pi : M_{i+1,i} = 0 \}. $$
  Then $B_i$ acts on $E^{[i,i+1]}_\pi$ by conjugation and $\cp$-conjugation,
  and the actions are the same. Also,
  $$ \mdeg E^\iip_\pi = (\A + z_{i+1} - z_i) \ \mdeg E_\pi.$$
\end{Lemma}

\begin{proof}
  Plainly the diagonal matrices act on $E^{[i,i+1]}_\pi$ with both actions
  the same, so it remains to check the action of $e^{i,i+1} \in Lie(B_i)$.
  If $N$ is the $\cp$-commutator of $e^{i,i+1}$ and $M$, then
  $$ N_{jk} \
  =\ [j=i]\ M_{i+1, k}\ [\cyc{i\leq i+1\leq k}]
  - M_{ji}\ [i+1 = k]\ [\cyc{j\leq i\leq i+1}] $$
  whereas the ordinary commutator doesn't have the cyclic requirements.
  But $[\cyc{i\leq i+1\leq k}]$ is $0$ only for $i=k$, and that term
  can safely be added since $M_{i+1, i}=0$ for $M\in E^{[i,i+1]}_\pi$.
  The other term is similar. So the $\cp$-commutator equals the
  ordinary commutator, hence the Lie algebra actions (and hence the
  Lie group actions) are the same.

  Since $B_i$ acts on $E$ by $\cp$-conjugation, and $B_i$ is connected,
  it acts on each component, such as $E_\pi$. 
  So it remains to check that $B_i$ preserves the subset $E^{[i,i+1]}_\pi$.
  For this we check the relevant matrix entry $N_{i+1,i}$ and see that
  is indeed zero.
  
  To see the claim about multidegrees, let $H$ be the space of
  matrices with vanishing $(i+1,i)$ entry (and vanishing diagonal). 
  By lemma \ref{lem:freeentries} and both parts of axiom (3') of multidegrees,
  $$ \mdeg E^\iip_\pi = (\A + z_{i+1} - z_{i}) \ \mdeg_H E^\iip_\pi
  = (\A + z_{i+1} - z_{i}) \ \mdeg E_\pi. $$
\end{proof}

We sequester some technicalities into a lemma whose proof we leave to
the reader.

\begin{Lemma}\label{lem:assoc}
%  \begin{eqnarray*}
%    [\cyc{i+1\leq l\leq m}]-[\cyc{i\leq l\leq m}] &=& \delta_{im}-\delta_{il}
%  \end{eqnarray*}
  Let $C,D \in \MNCdzero$. Then for all $j,m = 1\ldots N$,
  \begin{eqnarray*}
    ((e^{i+1,i} C)\cp D)_{jm}&=& (e^{i+1,i} (C\cp D))_{jm}
    + \delta_{i,m} \delta_{j,i+1} (CD)_{ii} \\
    (C\cp (D e^{i+1,i}))_{jm}&=& ((C\cp D)e^{i+1,i})_{jm} 
    + \delta_{i,m} \delta_{j,i+1} (CD)_{i+1,i+1} \\
    ((C e^{i+1,i})\cp D)_{jm}&=&  C_{j,i+1} D_{im}\ [\cyc{j\leq i\leq m}]\\
    (C\cp (e^{i+1,i} D))_{jm}&=&  C_{j,i+1} D_{im}\ [\cyc{j\leq i+1\leq m}] 
  \end{eqnarray*}
  If in addition $C_{i,i+1}=0$, then
  $$([e^{i+1,i},C])^{\cp 2}_{jm} = -\delta_{j,i+1}\delta_{i,m} (C^2)_{i,i+1}.$$
\end{Lemma}

We have reached the heart of the paper:

\begin{Proposition}\label{prop:recurrence}
  If the link pattern $\pi$ has no chord between $i$ and $i+1$, then
  \begin{equation}\label{eqn:transp}
    \mdeg E_{\pi} + \mdeg E_{f_i\cdot\pi}
  =-\ \frac{2\A + z_{i+1} - z_i}{\A + z_{i+1} - z_i}\ \der_i\ \mdeg E^\iip_{\pi}.
  \end{equation}
\end{Proposition}

\begin{proof} 
  The outline is as follows. We will apply lemma \ref{lem:divdiff} to
  the inclusion $f: E^\iip_\pi \to \MNC$, where $L_i$ acts on $\MNC$
  by ordinary conjugation. The image of $-\der_i f$ is ``almost'' inside $E$;
  by imposing one new equation (with weight $2\A+z_{i+1}-z_i$) the result $Z$
  is inside $E$. (There is an annoying technicality that we actually work
  not with $\Imgc{-\der_i f}$ but the larger and equally good 
  $Y := \Imgc{-\der_i f} + \complexes e^{i+1,i}$.) 
  Then we determine which components of $E$ are contained in $Z$:
  they are the $E_\pi$ and $E_{f_i\cdot\pi}$ components.

  We begin with the map $-\der_i f$.
  A generic element of $L_i$ can be written as $(1 + \lambda e^{i+1,i})b$
  for $b$ an element of $B_i$. We know by lemma \ref{lem:conjactions}
  that $B_i$ acting by conjugation preserves $E^\iip_\pi$. 
  To determine the closure of the image of $-\der_i f$, 
  it is therefore enough to look at the image of the dense open set
  $$ Q := 
  \{ (1 + \lambda e^{i+1,i}, M) : \lambda\in \complexes, M\in E^\iip_\pi \}. $$
%  $$ \{ (1 + \lambda e^{i+1,i}) M (1 + \lambda e^{i+1,i})^{-1} : 
%  M \in E^\iip_\pi \}. $$
  Let $N = M + \lambda (e^{i+1,i} M - M e^{i+1,i})$ 
  be in the image of $Q$. In particular, 
  \begin{equation}\label{eqn:Nentry}
     N_{i+1,i} = M_{i+1,i} - \lambda(M_{i+1,i+1} - M_{ii}) = 0. 
  \end{equation}
%  Let $M\in E^\iip_\pi$, $N = M + \lambda [e^{i+1,i},M]$,
%  so $N$ is a generic element of $\Img -\der_i f$. 
%  Note the algebra identity $A[A,B] + [A,B]A = A^2 B + B A^2$.
  Then writing $N^{\cp 2}$ for $N\cp N$, we have
  \begin{eqnarray*}
     N^{\cp 2} 
     &=& (M + \lambda [e^{i+1,i},M])^{\cp 2} \\
     &=& M \cp M + \lambda ([e^{i+1,i},M] \cp M + M \cp [e^{i+1,i},M])
     + \lambda^2 ([e^{i+1,i},M])^{\cp 2}  
  \end{eqnarray*}
  though the $\lambda^2$ term actually vanishes, thanks to
  the condition on $\pi$ and proposition \ref{prop:linearcondition}.
  Using lemma \ref{lem:assoc}, and $M\cp M = 0$, we get
  \begin{eqnarray*}
     (N^{\cp 2})_{jm}
     &=& \lambda ((e^{i+1,i} M) \cp M 
     - M \cp (M e^{i+1,i}) 
     + \lambda ([e^{i+1,i},M])^{\cp 2} 
     + M \cp (e^{i+1,i} M) 
     - (M e^{i+1,i}) \cp M)_{jm}   
     \\
     &=& \lambda(\delta_{i,m} \delta_{j,i+1} 
     ((M^2)_{ii} - (M^2)_{i+1,i+1} - \lambda (M^2)_{i,i+1}) \\
     &&\quad + ([j \leq i+1 \leq m] - [j\leq i\leq m])\ M_{j,i+1} M_{i,m} )
     \\
     &=& \lambda(\delta_{i,m} \delta_{j,i+1} 
     ((M^2)_{ii} - (M^2)_{i+1,i+1} - \lambda (M^2)_{i,i+1}) 
     + M_{j,i+1} M_{i,m} (\delta_{j, i+1} - \delta_{ m, i}))
     \\
     &=& \lambda \delta_{i,m} \delta_{j,i+1} 
     ((M^2)_{ii} - (M^2)_{i+1,i+1} - \lambda (M^2)_{i,i+1})
  \end{eqnarray*}
  where we used $M \in \MNCdzero$ to get from the third to the fourth line.
  So $N^{\cp 2}$ is zero away from $(N^{\cp 2})_{i+1,i}$.

  In particular, $N\in E$ if and only if $\lambda = 0$ 
  or $\lambda (M^2)_{i,i+1} = (M^2)_{ii} - (M^2)_{i+1,i+1}$. 
  In the latter case, 
  \begin{eqnarray*}
    (N^2)_{ii} &=& 
    (\exp(\lambda e^{i+1,i}) M \exp(-\lambda e^{i+1,i}))^2_{ii} \\
    &=& (\exp(\lambda e^{i+1,i}) M^2 \exp(-\lambda e^{i+1,i}))_{ii} \\
    &=& (M^2 + \lambda [e^{i+1,i},M^2])_{ii} \\
    &=& (M^2)_{ii} + \lambda (e^{i+1,i} M^2 - M^2 e^{i+1,i})_{ii} \\
    &=& (M^2)_{ii} - \lambda (M^2)_{i,i+1} \\
    &=& (M^2)_{i+1,i+1}
  \end{eqnarray*}
  Similarly $(N^2)_{i+1,i+1} = (M^2)_{ii}$, whereas $(N^2)_{jj} = (M^2)_{jj}$
  for all other $j$. Hence the only top components of $E$ that can appear in 
  the image of $-\der_i f$ are $E_\pi$ and $E_{f_i\cdot\pi}$. 
  
  We use this same calculation to prove that $-\der_i f$ is
  generically $1$:$1$, as it is enough to prove it on $Q$. Assume that
  $1$:$1$ness is violated:
  \begin{eqnarray*}
    &&-\der_i f(1 + \mu e^{i+1,i}, M) = -\der_i f(1 + \nu e^{i+1,i}, N) \\
    &\Longleftrightarrow&
    \exp(\mu e^{i+1,i}) M \exp(-\mu e^{i+1,i}) = 
    \exp(\nu e^{i+1,i}) N \exp(-\nu e^{i+1,i})  \\
    &\Longleftrightarrow&
    \exp((\mu-\nu)e^{i+1,i}) M \exp((\nu-\mu)e^{i+1,i}) = N \\
    &\Longrightarrow&
    \exp((\mu-\nu)e^{i+1,i}) M \exp((\nu-\mu)e^{i+1,i}) \in E_\pi
  \end{eqnarray*}
  As we saw above, there are only two possibilities for $\mu-\nu$ for
  which this left hand side is even in $E$, much less $E_\pi$. 
  If $M$ is in the open set 
  we used in theorem \ref{thm:topcomps} to define $E_\pi$, this
  $\exp((\mu-\nu) e^{i+1,i}) M \exp((\nu-\mu) e^{i+1,i})$ is not in $E_\pi$
  unless $\mu-\nu=0, M=N$. 

  We learn two things from $-\der_i f$ being generically $1$:$1$ :
  \begin{itemize}
  \item  $\dim \Imgc{-\der_i f} =
    \dim (L_i \times_{B_i} E^\iip_\pi)$, which in turn equals 
    $1 + \dim E^\iip_\pi = \dim E_\pi = \dim E$. (In fact this only used
    $-\der_i f$ being finite:$1$.)
  \item  by lemma \ref{lem:divdiff}, 
  $\mdeg \Imgc{-\der_i f} = -\der_i \mdeg E^\iip_\pi$.
  \end{itemize}

  Let $Y := \Imgc{-\der_i f} + \complexes e^{i+1,i}$. We note four properties
  of $Y$, the first two from equation (\ref{eqn:Nentry}):
  \begin{itemize}
  \item $\dim Y = \dim \Imgc{-\der_i f} + 1 = \dim E + 1$. 
  \item $ (\A+z_{i+1}-z_i)\, \mdeg Y = \mdeg \Imgc{-\der_i f} 
    = -\der_i\ \mdeg E^\iip_\pi \\ {\,}\hfill
    = -\der_i (\A+z_{i+1}-z_i)\, \mdeg E_\pi. $
  \item By lemma \ref{lem:freeentries} and calculation of $N^{\cp 2}$, 
    $M\in Y$ implies $(M\cp M)_{jm}=0$ unless $(j,m)=(i+1,i)$.
  \item $Y$ is irreducible, since $E_\pi$ hence $E^\iip_\pi$ hence
    $L_i \times_{B_i} E^\iip_\pi$ hence $\Imgc{-\der_i f}$ were.
  \end{itemize}

  Let $Z$ denote the intersection of $Y$ with the hypersurface
  $\{M : (M\cp M)_{i+1,i} = 0\}$, so $\dim Z \geq \dim Y - 1$. 
  Since $Y$ satisfied all but this one of $E$'s defining equations,
  $Z \subseteq E$, so $\dim Z \leq \dim E = \dim Y - 1$. 
  Hence $\dim Z = \dim E$, and
  $$ \mdeg Z = (2\A+z_{i+1}-z_i) \mdeg Y $$
 where $2\A+z_{i+1}-z_i$ is the $T$-weight of the equation of that hypersurface.
  Note that this is the right-hand side of the equation we seek.

  Since $Z\subseteq E$ and $\dim Z = \dim E$, the top-dimensional components
  of $Z$ are a selection of the top-dimensional components of $E$.
  Since $E$ is generically reduced on its top-dimensional components,
  $Z$ is too. 

  Plainly $E_\pi \subseteq Z$, since $-\der_i f$ restricted to 
  $B_i \times_{B_i} E^\iip_\pi$ already has image $E^\iip_\pi$. 
  We've already shown (by looking at the diagonal elements of the square)
  that the only other component of $E$ that could appear in $Z$ is
  $E_{f_i\cdot \pi}$.

  At this point we have two possibilities for the top components of $Z$:
  just $E_\pi$, or $E_\pi \cup E_{f_i\cdot \pi}$. Assume (for contradiction)
  the first. Then we have $\mdeg Z = \mdeg E_\pi$, so
  \begin{eqnarray*}
    \mdeg E^\iip_\pi &=& (\A+z_{i+1}-z_i) \mdeg E_\pi = (\A+z_{i+1}-z_i)\mdeg Z 
    \\
    &=& (2\A+z_{i+1}-z_i)(\A+z_{i+1}-z_i) \mdeg Y 
    = -(2\A+z_{i+1}-z_i) \der_i \mdeg E^\iip_\pi 
  \end{eqnarray*}
  Apply $\der_i$ to both sides:
  $$ \der_i \mdeg E^\iip_\pi = 2 \der_i \mdeg E^\iip_\pi $$
  so $\mdeg E^\iip_\pi = 0$, which is false by the positivity of multidegrees.

  Hence $Z$ has two top components, $E_\pi$ and $E_{f_i\cdot \pi}$, both
  generically reduced, so
  $$ \mdeg E_\pi + \mdeg E_{f_i\cdot\pi} = \mdeg Z 
=-\ \frac{2\A + z_{i+1} - z_i}{\A + z_{i+1} - z_i}\ \der_i\ \mdeg E^\iip_{\pi}.$$
\end{proof}

\begin{Theorem}\label{thm:main}
  The multidegree of $E_\pi$ is the homogenized component $\Psi_\pi$
%  of \cite{DFZJ}
  for all link patterns $\pi$:
  \begin{equation}    \label{eqn:mainresult}
    \mdeg E_\pi|_{\A=1}=\Psi_\pi
\end{equation}
\end{Theorem}

\begin{proof}
  Setting $\A=1$ in equations~(\ref{eqn:pizero}) and (\ref{eqn:transp}), one
  recovers the equations~(\ref{eqn:pizeroa}) and (\ref{eqn:schub}).
  As explained in \cite{DFZJ} these characterize the $\Psi_\pi$ uniquely,
  hence the equality (\ref{eqn:mainresult}).
\end{proof}

This has a corollary, for which it would be interesting to have a
geometric proof along the lines of proposition \ref{prop:recurrence}.
%{\bf check for $N$ odd}

\begin{Corollary}\label{cor:smallarch}
  If $\pi(i)=i+1$, then
  $$-\der_i\ (\A+z_i-z_{i+1}) \mdeg E_\pi^\iip
  = -2\A \sum_\rho \der_i\ \mdeg E_\rho^\iip $$
  where the sum is taken over those $\rho \neq \pi$ such that 
  $e_i\cdot \rho = \pi$ and the $i,i+1$ strands cross.
\end{Corollary}

\begin{proof}
  Let $u=z_i-z_{i+1}$.
  The $\pi$-component of the equations (\ref{eqn:Rcheck})
  $$ a(u) \Psi_\pi + b(u) \Psi_{f_i\cdot \pi}
  + c(u) \sum_{\rho:\ e_i\cdot \rho = \pi} \Psi_\rho 
  = \tau_i \Psi_\pi $$
  looks different depending on whether $\pi(i) = i+1$ or $\pi(i)\neq i+1$.  
  The $\pi(i)\neq i+1$ equations are simpler because the summation
  term vanishes, and can be rewritten as the equations (\ref{eqn:schub}).

  If however $\pi(i)=i+1$, we rewrite as
  $$ \big(a(u)+b(u)+c(u)\big) \Psi_\pi 
  + c(u) \sum_{\rho:\ e_i\cdot \rho = \pi, \rho\neq \pi} \Psi_\rho 
  = \tau_i \Psi_\pi. $$
  Substituting in the formula (\ref{eqn:Rcheckform}) for $\check R$, 
  this becomes
  $$  \Psi_\pi   + \frac{2u}{(2-u)(1+u)} 
  \sum_{\rho:\ e_i\cdot \rho = \pi, \rho\neq \pi} \Psi_\rho 
  = \tau_i \Psi_\pi $$
   Hence
  $$ -\der_i \Psi_\pi = \frac{1}{u} (\tau_i \Psi_\pi - \Psi_\pi) =
  \frac{2}{(2-u)(1+u)} \sum_{\rho:\ e_i\cdot \rho = \pi, \rho\neq \pi} 
  \Psi_\rho. $$
  The $\rho$-terms in the summation can be grouped into pairs 
  $\{\rho,f_i\cdot \rho\}$. Since $\rho\neq \pi$, no $\rho=f_i\cdot \rho$,
  and $\rho(i)\neq i+1$. We can pick a preferred element of each pair
  by asking that the $(i,\rho(i))$ chord cross the $(i+1,\rho(i+1))$.
  Then the equation becomes
  \begin{eqnarray*}
    -\der_i \Psi_\pi 
    &=& \frac{2}{(2-u)(1+u)}\sum_{\rho} (\Psi_\rho + \Psi_{f_i\cdot \rho})\\
    &=& \frac{2}{(2-u)(1+u)}\sum_{\rho}-\frac{2-u}{1-u}\der_i (1-u)\Psi_\rho\\
    &=& \frac{-2}{(1+u)(1-u)}\sum_{\rho}\der_i (1-u)\Psi_\rho
  \end{eqnarray*}
  where the summation is over 
  $\rho\neq \pi$, $e_i\cdot \rho = \pi$, 
  and the $i,i+1$ chords of $\rho$ cross. So
  $$ -\der_i (1+u)(1-u)\Psi_\pi = -(1+u)(1-u) \der_i \Psi_\pi = 
  -2\sum_{\rho} \der_i (1-u)\Psi_\rho $$
  Using theorem \ref{thm:main} we obtain the desired formula.
\end{proof}

We will give a direct geometric derivation of this result in \cite{DFKZJ}.

\section{The permutation sector and the upper-upper scheme}\label{sec:perm}

In this section we work again in the $(R,L)$ coordinate system on $E$.

Define the {\dfn permutation subspace} $\MNC_P$ to be the
subspace of $\{(R,L)\}$ in which the upper triangular matrix $R$ 
is supported in the northeast rectangle:
$$ R_{ij} = 0 \quad \hbox{unless} \quad 
i\leq n,\ j\geq n+1 $$
(recall that $N=2n+r$, $r=0$, $1$). 
It is easy to check that $\MNC_P$ is invariant
under $\cp$-conjugation by $U$.

Let $X$ denote that northeast rectangle (or square, if $N$ even), 
so $X$ is an $n \times (n+r)$ matrix with
$$ X_{ij} = R_{i,j+n}. $$ 
Similarly, let $Y$ denote the transposed rectangle in $L$, 
so $Y$ is an $(n+r) \times n$ matrix with
$$ Y_{ij} = L_{i+n,j}. $$ 
Put together, 
\begin{equation}
  \label{eq:MtoXY}
   M =
\begin{pmatrix}
  *\backslash 0
  &  & X \\
  Y &  &
  *\backslash 0
\end{pmatrix}
\qquad
\hbox{where each $*\backslash 0$ is strictly lower triangular.}
\end{equation}

Define the {\dfn permutation sector} $E_P\subseteq E$ to be the 
intersection $E \cap \MNC_P$. Then (as in lemma \ref{lem:freeentries})
the conditions on $R$ and $L$ are in fact only conditions on $X$ and $Y$:
$$ \forall (R,L)\in \MNC_P, \qquad (R,L) \in E_P \Longleftrightarrow
XY, YX \hbox{ are upper triangular square matrices}. $$
(Note that if $N$ is odd, then $YX$ is one size larger than $XY$.)
In the case $N$ even, this ``upper-upper scheme'' $E_P$
was introduced in \cite{Kn}, and most of the next theorem proven.
The case $N$ odd was considered in \cite{dGN}.

Note that since $E_P$ lives inside the linear subspace $\MNC_P$,
its multidegree and that of its components are divisible by
$$ \mdeg \MNC_P = \prod_{1\leq i\leq j\leq n} (\A + z_i - z_j)
   \prod_{n+1\leq i\leq j\leq N} (\A + z_i - z_j).
$$
With these factors divided out, we recover the multidegrees
relative to $\MNC_P$.

Let $P$ denote the set of link patterns $\pi$ such that
$\forall i=1,\ldots,n$, $\pi(i)>n$.
For $i>n$, this forces $\pi(i)\leq n$ or ($N$ odd) $\pi(i)=i$.
For $N$ even (the case considered in \cite{DFZJ}),
such $\pi$ correspond in an obvious way to permutations of
$\{1,\ldots,n\}$.

\begin{Theorem}
  The permutation sector $E_P$ is a complete intersection,
  hence has multidegree
  $$
  \mdeg E_P =
  \prod_{1\leq i\leq j\leq n} (\A + z_i - z_j) (2\A + z_j - z_i)
  \prod_{n+1 \leq i\leq j\leq N} (\A + z_i - z_j) (2\A + z_j - z_i)
  $$
  as a subscheme of $\MNCdzero$.

  Moreover, $E_P = \cup_{\pi \in P} E_\pi$, and in particular is reduced.
\end{Theorem}

\begin{proof}
  A complete intersection, by definition, is a scheme $C$ whose codimension 
  equals the number of defining equations. It is enough to check that
  $\codim C$ is at least this number of equations, as
  the inequality then implies the equality.
  There are
  $$ \bigg({\lfloor N/2\rfloor \choose 2} + {\lceil N/2\rceil \choose 2}\bigg) 
  + {\lfloor N/2\rfloor \choose 2} + {\lceil N/2\rceil \choose 2}
  $$
  equations, 
  for the vanishing of the two parts of $R$ outside $X$, the strict lower
  triangle of $XY$, and the strict lower triangle of $YX$.
  
  Since $E_P$ is a subscheme of $E$, its codimension is at least that of $E$,
  namely $(N^2 - N) - \lfloor N^2/2 \rfloor = \lceil N^2/2 \rceil - N$. 
  (We are computing codimension relative 
  to the $(N^2-N)$-dimensional vector space $\MNCdzero$,
  and using theorem \ref{thm:reduced}.)
  If $N=2n$, then the number of defining equations and codimension are
  $4{n\choose 2} \leq 2n^2 - 2n$. If $N=2n+1$, then these two numbers
  are $2 {n\choose 2} + 2{n+1 \choose 2} \leq 2n^2+2n+1 - (2n+1)$.
  In either case we get the desired inequality (with, of course, equality).

  The multidegree of a complete intersection of $T$-invariant hypersurfaces
  is the product of the weights of the defining equations. This gives
  the stated multidegree for $E_P$.

  We make use of two properties of complete intersections: they are
  equidimensional, and more specifically Cohen--Macaulay.
  Since $E_P$ is equidimensional, its support is a
  union of components of $E$. Since $E$ is generically reduced, so is $E_P$.
  Since $E_P$ is Cohen--Macaulay and generically reduced, it is reduced
  (this was already proven in \cite{Kn} in the case $N$ even, via
  the same argument). So scheme-theoretically it is the union of some
  components of $E$. % (which were reduced).
  
  Finally it remains to determine which components of $E$ lie in
  $E_P$.  By the proof of theorem \ref{thm:diagpaired}, we know that
  for generic diagonal $t$, we have $\Mpi t \in E_\pi$, $\Mpi t\notin
  E_\rho$ for $\rho\neq \pi$.  Since $E$ and $\MNC_P$ are $U$-invariant,
  so is $E_P$, hence $U\cdot \{\Mpi t\} \subseteq E$ for $\pi\in P$.
  Hence by proposition \ref{prop:openset}, $E_\pi \subseteq E_P$ if
  and only if $\Mpi t \in E_P$. The vanishing conditions on $E_P$ are
  then equivalent to $\pi \in P$.
\end{proof}

By the additivity of multidegrees, and theorem \ref{thm:main}, we have the

\begin{Corollary}\cite[for $N$ even]{DFZJ}
  $$ \sum_{\pi \in P} \Psi_\pi = 
  \prod_{1\leq i\leq j\leq n} (1 + z_i - z_j) (2 + z_j - z_i)
  \prod_{n+1\leq i\leq j\leq N} (1 + z_i - z_j) (2 + z_j - z_i)
  $$
\end{Corollary}
The $N$ unrestricted, $z_i \equiv 0$ case was conjectured in \cite{dGN}.

Finally, we prove the original observation of \cite{dGN}, i.e.\ 
that the component $\Psi_\pi$ for $\pi(i)=2n+1-i$ provides the
degree of the commuting scheme $C_n=\{ (X,Y)\in \MnC : XY=YX\}$.
We strengthen this to a computation of the multidegree (for a new torus, 
as not all of $T$ acts on $C_n$).

\begin{Theorem}\label{thm:comvar}
  Let $S$ be the product of $\complexes^\times$
  and the diagonal matrices in $\MnC$. Then $S$ acts on $C_n$ by 
  $$ (\alpha, D) \cdot (X,Y) := (\alpha D X D^{-1}, \alpha D Y D^{-1}). $$

  Let $N=2n$. 
  Let $\Phi_n = \mdeg \MNC_P = \prod_{1\leq i\leq j\leq n} (\A + z_i - z_j)
  \prod_{n+1\leq i\leq j\leq N} (\A + z_i - z_j)$.
  Define $\pi_n$ by $\pi_n(i):=N+1-i$. 

  Then the $S$-multidegree of the commuting scheme inside $\MnC\times\MnC$
  satisfies
  $$  S\!-\!\mdeg_{\MnC\times \MnC} C_n 
  = (\mdeg E_{\pi_n} / \Phi_n)|_{z_i\equiv z_{n+i}}, $$
  where we 
  denote the standard basis for $S$'s weight lattice by $(\A,z_1,\ldots,z_n)$.

%  so the weights on $\MnC\times\MnC$ are $\{\A+z_i-z_j : i,j=1\ldots n\}$ 
%  each with multiplicity $2$.
  
  Setting $\A=1, z_i \equiv 0$ on both sides, we get $\deg C_n=\deg E_{\pi_n}$.

  Let $\Delta_n := (\mdeg E_{\pi_n}/\Phi_n)|_{z_{n+1} = \cdots = z_{2n} = 0}$,
  so $\deg C_n = \Delta_n|_{\A=1,z_{1} = \cdots = z_{n} = 0}$.
  This $\Delta_n$ can be calculated as
  \begin{eqnarray*}
    \Delta_n
    &=& \A^n\ 
      \theta_1\ \theta_2 \theta_1\ \cdots\ \theta_i
      \theta_{i-1}\cdots \theta_2 \theta_1 \ \cdots\ 
      \theta_{n-1}\cdots\theta_2\theta_1 
      \prod_{i=1}^n (\A+z_i)^{i-1}(\A-z_i)^{n-i}  \\
    &=& \A^n\ \theta_1\theta_2 \cdots\theta_{n-1}\ 
      \theta_1 \theta_2 \cdots \theta_{n-2}\ 
      \cdots\ 
      \theta_1 \theta_2 \ \theta_1
      \prod_{i=1}^n (\A+z_i)^{i-1}(\A-z_i)^{n-i}  \\
    &=& \A(\theta_1 \cdots \theta_{n-1}) (\A+z_n)^{n-1} \prod_{i=1}^{n-1} (\A-z_i)
    \ \Delta_{n-1}
  \end{eqnarray*}
  where $\theta_i=-2\A\der_i-\tau_i$.
\end{Theorem}

\begin{proof}
  In \cite{Kn}, the equations of the commuting scheme $C$ are 
  $S$-equivariantly degenerated to those of
  $$  F_+ := \{ (X,Y) : XY, YX \mbox{ upper triangular, } 
  \diag(XY) = \hbox{ reverse of }\diag(YX) \}. $$
  (The details of this family are unimportant here.)
  While this degeneration is conjectured in \cite{Kn} to be flat, 
  this is not proven. 
  So a priori one only knows that the actual flat limit $F$ of the
  commuting scheme is contained inside $F_+$. 
  
  This upper bound $F_+$ is contained in the upper-upper scheme, 
  and it is easy to check that it contains one entire component $F_-$
  (corresponding to the reversal permutation)
  and only lower-dimensional parts of other components.
  Since the upper-upper scheme is generically reduced, 
  $F_+$ is generically reduced along $F_-$. 
  Putting these two facts together and applying axiom (2) of multidegrees,
  we see that $F_+$ and $F_-$ have the same $S$-multidegree.
  
  Since the upper-upper scheme has the same dimension, $n^2+n$, as the
  commuting scheme, $n^2+n \geq \dim F_+ \geq \dim F = n^2+n$.
  
  Since $F$ is a degeneration of the (irreducible) commuting scheme, 
  it is set-theoretically equidimensional. The only component of $F_+$ of the 
  right dimension is $F_-$, so $F \supseteq F_-$ and they are equal as sets.
  Since $F_+$ and $F_-$
  have the same $S$-multidegree, $F$ trapped between them has the 
  same $S$-multidegree as both. 
  (If $F = F_-$, this is enough to prove that $C_n$ is reduced, which is
  still unknown. In \cite{Kn} it is further conjectured that 
  $F_+ = F = F_-$.) This is also the $S$-multidegree of $C_n$, since $C_n$
  degenerates to $F$. 

  Embed $\MnC\times\MnC$ into $\MNC$ as in equation (\ref{eq:MtoXY}).
  Our $(N+1)$-dimensional torus $T$ acts on $\MnC\times \MnC$, preserving
  $F_+$ and $F_-$, by
  $$ (\alpha,D_1,D_2) \cdot (X,Y) :=
  (\alpha D_1 X D_2^{-1}, \alpha D_2 Y D_1^{-1}) $$
  where
  $D_1 :=\diag(\zeta_1,\ldots,\zeta_n), 
  D_2  :=\diag(\zeta_{n+1},\ldots,\zeta_{2n})$.
  This extends the action of $S$, which corresponds to the subtorus $D_1=D_2$.
  Hence
  $$   \mdeg_{\MnC\times \MnC} F_-|_{z_i\equiv z_{n+i}} 
  = S\!-\!\mdeg_{\MnC\times\MnC} F_-. $$

  The component of $E_P$ corresponding to the component $F_-$ of the
  upper-upper scheme is $E_\pi$ where $\pi(i) = 2n+1-i$. Hence
  $$ \mdeg E_{\pi_n} = \Phi_n\ \mdeg_{\MnC\times\MnC} F_-. $$
  Combining these equations,
%  $S\!-\!\mdeg_{\MnC\times\MnC} C_n = S\!-\!\mdeg_{\MnC\times\MnC} F_-$,
  we get
  \begin{eqnarray*}
    S\!-\!\mdeg_{\MnC\times\MnC} C_n &=& S\!-\!\mdeg_{\MnC\times\MnC} F_- \\
    &=& \mdeg_{\MnC\times \MnC} F_-|_{z_i\equiv z_{n+i}} \\
    &=& (\mdeg E_{\pi_n}/\Phi_n)|_{z_i\equiv z_{n+i}}.
  \end{eqnarray*}

  To compute $\mdeg E_{\pi_n}/\Phi_n$,
  we apply equation~(\ref{eqn:schub}) in order to
  produce the pattern $\pi$ out of $\pi_0$ using transpositions
  $f_i$ with $i$ taking values in $\{1,\ldots,n-1\}$ only; %so that
  the variables $z_i$, $i=n+1,\ldots,2n$ do not appear in the
  divided difference operators used and can be set to zero from the start. 
  Hence we can use the same calculation to compute $\Delta_n$. (In fact
  $\Delta_n$ can also be interpreted as a multidegree, for the
  subtorus $D_2=1$, but this won't be relevant for us.)
  Finally, one can pull out a factor $\prod_{1\le i<j\le n} (\A+z_i-z_j)$, 
  turning the divided difference operator of equation (\ref{eqn:schub}) 
  into $\theta_i$ (cf. \cite[equation (4.20)]{DFZJ}).
  The first two formulae (which only differ in the order of the
  $\{\theta_i\}$) follow. As explained in the introduction, the first
  is more practical.
  
  If we expand out the recursive formula given for $\Delta_n$, we get
  $$
  \Delta_n  = 
  \A(\theta_1 \cdots \theta_{n-1}) (\A+z_n)^{n-1} \prod_{i=1}^{n-1} (\A-z_i)
  \ \ 
  \A(\theta_1 \cdots \theta_{n-2}) (\A+z_{n-1})^{n-2} \prod_{i=1}^{n-2}(\A-z_i)
  \ \ \cdots
  $$
  To get from there to the second closed form, notice first that
  $\theta_i(pq) = p\theta_i(q)$ if $p$ is symmetric in $\{z_i,z_{i+1}\}$.
  Hence we can pull each of the multiplicative factors
  $(\A+z_m)^{m-1} \prod_{i=1}^{m-1} (\A-z_i)$ to the right, past each
  $\theta_j$, since they only encounter $\theta_j$ for $j<m-1$.
  Similarly, pull the isolated factors of $\A$ left.
  $$
  \Delta_n  = \A^n
  (\theta_1 \cdots \theta_{n-1}) (\theta_1 \cdots \theta_{n-2}) \cdots
  (\A+z_n)^{n-1} \bigg( \prod_{i=1}^{n-1} (\A-z_i)\bigg)
  (\A+z_{n-1})^{n-2} \bigg( \prod_{i=1}^{n-2} (\A-z_i) \bigg)
  \cdots
  $$
  Reordering the multiplicative factors, we get the second closed formula.
\end{proof}

\section{From $N$ to $N-2$ as a geometric vertex decomposition}

In this section we give a geometric interpretation
of theorem 4 from \cite{DFZJ}.
Its proof is based on corollary 2.5 of \cite{KMY}, which reads

\newcommand\Sym{{\rm Sym}}
\begin{Corollary*}
  Let $X \subseteq H\times L$, where $H$ has coordinates $x_1,\ldots,x_n$
  and $L$ has coordinate $y$. Assume that $H$, $L$ are representations
  of a torus $T$, and $X$ is a $T$-invariant subvariety. 
  Let $w \in T^* = \Sym^1(T^*)$ 
  be the weight of $T$ on $L$, and $S\leq T$ the stabilizer of $L$, so
  the map $\Sym(T^*)\to \Sym(S^*)$ takes $p\mapsto p|_{w=0}$.

  Let the ideal $I$ defining $X$ be generated
  by $\{y^{d_i} q_i + r_i\}_{i=1}^m$, where $y^{d_i} q_i$ is sum of the
  terms in $y^{d_i} q_i + r_i$ divisible by the highest power of $y$.
  Let $J = \< q_i\ |\ i=1,\ldots,m \>$. 
  Let $\Theta \subseteq H$ be the corresponding scheme.
  If we know that 
  \begin{itemize}
  \item $\Theta$ has only one component of dimension $\dim X - 1$
  \item that component is generically reduced
  \item $X$ is not contained in a union of finitely many translates of $H$
  \end{itemize}
  then
  $$ (\mdeg_{H\times L} X)|_{w=0} = (\mdeg_H \Theta)|_{w=0}. $$
\end{Corollary*}

Thanks to theorem \ref{thm:main} relating $\mdeg E_\pi$ to $\Psi_\pi$,
the following is exactly theorem 4 from \cite{DFZJ}. We explain after
the theorem what new insight is available from a geometric proof.

\begin{Theorem}
  Let $\pi$ be a link pattern on $1\ldots N$, and $i$ such that
  $\pi(i)=i+1$. We can associate to it a smaller link pattern $\pi'$ on 
  $1\ldots i-1,i+2,\ldots N$.

  Then if we specialize $\mdeg E_\pi$ at $z_{i+1}=z_i+\A$, we get
  $$ \mdeg E_\pi\big |_{z_{i+1}=z_i+\A} =
 \bigg(   \prod_{k \neq i,i+1} (\A+z_{i+1}-z_k)(\A+z_k-z_i)   \bigg)
  \  \mdeg_{M_{N-2}(\complexes)_{\Delta=0}}
  E_{\pi'}(\ldots,z_{i-1}, z_{i+2},\ldots).$$ 
\end{Theorem}

\begin{proof}
  We use the notation of the corollary throughout.
  Let $H\times L = \MNCdzero$,
  where $H = \{ M \in \MNCdzero: M_{i,i+1} = 0\}$ and
  $L %= \{ M \in \MNC: M_{kl} \neq 0 \implies (k,l)=(i,i+1)\}
  =\complexes e^{i,i+1}$. Then
  $w = \A+z_i-z_{i+1}$.

  The equations
  $(M\cp M)_{ab}=\sum_{j:\ \cyc{a\le j\le b}} M_{aj} M_{jb}$
  involve the entry $M_{i,i+1}$ if and only if $a=i$ or $b=i+1$, 
  and their $y^{d_i} q_i$ terms are 
  $ M_{i,i+1} M_{i+1,b} = 0$ for $b \neq i$,
  $ M_{ai} M_{i,i+1} = 0$ for $a \neq i+1$.
  This gives us some linear equations on $\Theta$:
  $$  M_{i+1,b} = 0\quad \hbox{for }b\neq i,
  \qquad M_{ai} = 0\quad \hbox{for }a\neq i+1. $$
  In particular, there are no equations on $\Theta$ involving the
  entries $M_{i,*},M_{*,i+1}$ except the diagonal entries.

  Using the equations from theorem \ref{thm:compeqns},
  we can determine enough of the remaining equations on $\Theta$ 
  to compute its multidegree:
  given $M\in \Theta$, if we let $M'$ be $M$ with its $i$th and $(i+1)$st 
  rows and columns erased, then $M'$ satisfies 
  the equations from theorem \ref{thm:compeqns} on $E_{\pi'}$. 
  Applying axiom (3')
  of multidegrees, we get a linear factor for each vanishing
  $\{ M_{i+1,b} , M_{ai} \}$, and the relation
  $$ \mdeg \Theta = 
  \bigg(   \prod_{k \neq i,i+1} (\A+z_{i+1}-z_k)(\A+z_k-z_i)   \bigg)
  \mdeg_{M_{N-2}(\complexes)_{\Delta=0}} 
  E_{\pi'}(1,\ldots,z_{i-1}, z_{i+2},\ldots, z_N).$$ 
  Then apply the corollary.
\end{proof}

Put another way, the difference between $\mdeg E_\pi$ and the right-hand
side of this equation is a multiple of $\A+z_i-z_{i+1}$. 
With \cite[theorem 2.4]{KMY}, we can give an interpretation of that multiple:
it is the multidegree of the projection of $E_\pi$ to the subspace
$H = \{ M \in \MNCdzero: M_{i,i+1} = 0\}$.

\section{The flat limit $D_0$ of the nilpotent orbit $\{M: M^2=0\}$}
\label{sec:flatlimit}

In this section we elucidate the precise relation between $E$ and the
scheme $D_1 := \{M\in \MNC : M^2 = 0\}$.
We begin with some results about $D_1$.

\begin{Lemma}
  $D_1$ is irreducible.
  For $N$ even, $D_1$ is generically reduced; for $N$ odd,
  it contains the underlying reduced scheme with multiplicity $2$.
\end{Lemma}

\begin{proof}
  The fact that $D_1$ is irreducible follows from Jordan canonical form.
  To check if $D_1$ is generically reduced we consider the point $M$
  with the following block structure:
  $$ M=\begin{pmatrix}0_{(n+r)\times n}&0_{(n+r)\times(n+r)}\\ 
    {\bf 1}_{n\times n}&0_{n\times(n+r)}\end{pmatrix}$$
  The Zariski tangent space is the kernel of $L\mapsto ML+LM$. For $N$
  even, this leads to the set of equations
  $$ L_{ij}=-L_{i+n,j+n}\quad {\rm and}\quad L_{i,j+n}=0 \quad
  i,j=1\ldots n$$
  hence the correct codimension of $2n^2$, which implies the generic
  reducedness of $E_N$.  For $N$ odd, we find this time
  $$L_{ij}=-L_{i+n+1,j+n+1}\quad i,j=1\ldots n\quad{\rm and}\quad
  L_{ij}=0 \quad i=1\ldots n+1, j=n+1\ldots N, (i,j)\ne (n+1,n+1)$$
  hence a codimension of $n^2+(n+1)^2-1=2n(n+1)$ which is one less
  than the codimension of $E_N$.  Note however that adding the extra
  equation $\Tr M=0$ $\Rightarrow$ $\Tr L=0$ increases the codimension
  by $1$ and makes $E_N$ generically reduced. We now show that
  generically $(\Tr M)^2=0$, thus the multiplicity is $2$.
  
  This requires a bit more work, since we must go back to a generic
  $M$. We consider the following matrix $P$ given by
  $$ P_{ij}=
  \begin{cases}
    \delta_{ij}& i\le n+1\\
    M_{i-n,j}&i>n+1.
  \end{cases}
  $$
  Generically, $\det P$ is non-zero on $D_1$. 
  (Otherwise, it would be identically zero since $D_1$ is irreducible,
  but it is easy to construct an $M\in D_1$ for which $\det P\neq 0$.)
  We therefore allow ourselves to invert $\det P$, and in particular
  to use the inverse matrix $P^{-1}$.  Thanks to $M^2=0$, $PMP^{-1}$ has
  a certain block structure which can be summarized as follows:
  $$ (PMP^{-1})_{ij}=
  \begin{cases}
    u_j& i=1\\
    \delta_{i-n,j}& 1<i\le n+1\\
    0& i>n+1\\
  \end{cases}
  $$
  where the $u_j$ are some polynomials of the $M_{ij}$ and of $\det
  P^{-1}$ whose explicit form is not needed.  Note that this is not
  quite the block structure of the $M$ chosen in the beginning of the
  proof.  This is because we have ``missed'' the fact that in odd
  dimension $\dim \Img M$ is generically one less than $\dim \Ker M$.
  We now write $(PMP^{-1})^2=PM^2 P^{-1}=0$ and indeed find the extra
  condition that $u_1^2=0$.  But $\Tr M=\Tr (PMP^{-1})=u_1$, hence
  generically $(\Tr M)^2=0$.
\end{proof}

In fact the radical of $D_1$'s ideal is generated by the entries of
$M^2$ and $M$'s characteristic polynomial \cite{St,We}.

\begin{Proposition}\label{prop:degD1}
  The multidegree of the scheme $D_1$ is
  \begin{eqnarray*}
    \mdeg_{\MNC} D_1 
    &=& 2^r \prod_{i, j}\ (\A+z_i-z_j) 
    \sum_{S \subseteq \{1,\ldots,N\} \atop |S|=n}
    \prod_{s\in S, \bar s\notin S} ((\A+z_s-z_{\bar s})(z_{\bar s}-z_s))^{-1}\\
    &=&   2^{n+r} \A^N\prod_{i<j} {(\A-z_i+z_j)(\A-z_j+z_i)\over (z_i-z_j)}
  {\rm Pf}\left({z_i-z_j\over (\A+z_i-z_j)(\A+z_j-z_i)}\right)_{1\le i,j\le N}.
  \end{eqnarray*}
  Moreover, the sequence $\{ \mdeg D_1^N \}$ (where the size $N$ of the
  matrices now varies) is characterized by the properties
  \begin{itemize}
  \item $\mdeg {D_1^0} = 1$, $\mdeg {D_1^1} = 2\A$
  \item $\mdeg D_1^N$ is a symmetric polynomial in $z_1,\ldots,z_N$
  \item 
    $$\mdeg_\MNC D_1^N\bigg|_{z_2=z_1+\A} 
    = 2 \A^2\left(\prod_{k=3}^N (\A+z_2-z_k)(\A+z_k-z_1)\right)
    \mdeg_{M_{N-2}(\complexes)} D_1^{N-2}\big(z_3,\ldots,z_N\big) $$
  \end{itemize}
\end{Proposition}

\newcommand\Gr{{\rm Gr}}
\newcommand\Hom{{\rm Hom}}
\newcommand\e[1]{\exp\left({#1}\right)}

\begin{proof}
  Let $Q = \{(V \in \Gr_n(\complexes^N), H \in \Hom(\complexes^N/V,V)) \}$
  be the vector bundle over the Grassmannian of $n$-planes in 
  $\complexes^N$, where the fiber over $V$ is the linear space of maps
  from $\complexes^N/V\to V$. (In fact $Q$ is isomorphic to the
  cotangent bundle.) Then there is a generically $1:1$ map
  \begin{eqnarray*}
     \beta : Q &\to& D_1 \\
     (V,H) &\mapsto& 
     (\complexes^N \onto \complexes^N/V \stackrel H \To
     V \into \complexes^N)
  \end{eqnarray*}
  which is equivariant with respect to an action of 
  $\complexes^\times \times GL_N(\complexes)$: let $\complexes^\times$
  rescale the fibers of the bundle and rescale $D_1$, and $GL_N(\complexes)$
  act in the obvious ways. In particular $\beta$ is equivariant for the
  action of our torus $T$. 
  %It is onto the {\em set} $D_1$, but does not
  %see the multiplicity $2$ for odd $N$.

  The $T$-fixed points on $Q$ are of the form $(V,\vec 0)$ where $V$ is
  an $n$-dimensional coordinate subspace $\complexes^S$, 
  using the coordinates $S\subseteq \{1,\ldots,N\}$.
  The tangent space
  $T_{(V,\vec 0)} Q$ is isomorphic to 
  $\Hom(V,\complexes^N/V) \oplus \Hom(\complexes^N/V,V)$,
  where the rescaling circle only acts on the second factor. The
  weights are $\{z_j-z_i\},\{\A+z_i-z_j\}$ where $i\in S,j\notin S$.

  Via the same sort of equivariant localization arguments as in 
  lemma \ref{lem:divdiff}, we obtain the formula
  $$ 1 =  \sum_S [(\complexes^S,0)]
  \prod_{i\in S, j\notin S} ((\A+z_i-z_j)(z_j-z_i))^{-1} $$
  as a formula in (a localization of) $H^*_T(Q)$, where $[(\complexes^S,0)]$
  is the class of the point $(\complexes^S,0) \in Q$. Pushing that into
  $\MNC$ using $\beta_*$, the class of each point maps
  to the class of $\{$the zero matrix$\}$, which is the product of the weights 
  on $\MNC$. Including the factor $2^r$ for the scheme structure,
  we get the desired formula and a close equivalent:
  \begin{eqnarray*}
    \mdeg_{\MNC} D_1 
    &=& 2^r \prod_{i, j}\ (\A+z_i-z_j)
    \sum_{S \subseteq \{1,\ldots,N\} \atop |S|=n}
    \prod_{i\in S, j\notin S} ((\A+z_i-z_j)(z_j-z_i))^{-1} \\
    &=& 2^r 
    \sum_{S \subseteq \{1,\ldots,N\} \atop |S|=n}
    \prod_{i,j \in S} (\A+z_i-z_j)
    \prod_{i,j \notin S} (\A+z_i-z_j)
    \prod_{i\in S,j \notin S} \frac{\A+z_j-z_i}{z_j-z_i}
  \end{eqnarray*}
  
  The base cases are obvious, and the symmetry follows from the
  $GL_N(\complexes)$ and hence $S_N$ action.
  We will see the recurrence relation from the second version of the
  formula above. If $\A+z_1-z_2=0$, the only nonzero terms have
  $S\ni 1, S\not\ni 2$, so we can separate out the factors involving $1,2$
  and rewrite 
  \begin{eqnarray*}
    \prod_{i,j \in S} (\A+z_i-z_j)&=& 
    \A \prod_{i\in S\setminus 1} (z_2-z_i)(\A+z_i-z_1)
    \prod_{i,j \in S\setminus 1} (\A+z_i-z_j)
    \\
    \prod_{i,j \notin S} (\A+z_i-z_j) &=&
    \A \prod_{j\notin S\cup 2} (\A+z_2-z_j)(z_j-z_1)
    \prod_{i,j \notin S\cup 2} (\A+z_i-z_j) 
    \\
    \prod_{i\in S,j \notin S} \frac{\A+z_j-z_i}{z_j-z_i} &=&
    2 \prod_{i\in S\setminus 1} \frac{\A+z_2-z_i}{z_2-z_i}
    \prod_{j\notin S\cup 2} \frac{\A+z_j-z_1}{z_j-z_1}
    \prod_{i\in S\setminus 1, j\notin S\cup 2} \frac{\A+z_j-z_i}{z_j-z_i}
  \end{eqnarray*}
  giving a total product of 
  \begin{eqnarray*}
    \mdeg_{\MNC} D_1^N %\bigg|_{z_2=z_1+\A} 
    &=& 2 \A^2 
    \sum_{S \subseteq \{1,3,\ldots,N\} \atop S\ni 1, |S|=n}
    \prod_{i\in S\setminus 1} (\A+z_i-z_1)(\A+z_2-z_i)
    \prod_{j\notin S\cup 2} (\A+z_2-z_j) (\A+z_j-z_1) \\
    && 2^r 
    \prod_{i,j \in S\setminus 1} (\A+z_i-z_j)
    \prod_{i,j \notin S\cup 2} (\A+z_i-z_j) 
    \prod_{i\in S\setminus 1, j\notin S\cup 2} \frac{\A+z_j-z_i}{z_j-z_i} \\
    &=& 2\A^2 \left(\prod_{i\neq 1,2} (\A+z_i-z_1)(\A+z_2-z_i)\right)
    \ \mdeg_{\MNC} D_1^{N-2}(z_3,\ldots,z_N)
  \end{eqnarray*}
  at $z_2=z_1+\A$, as desired.
  
  In \cite[theorem 5]{DFZJ} it was shown that the symmetry, base case,
  and recurrence relation are enough to determine $\sum_\pi \Psi_\pi$,
  and to derive a Pfaffian formula. The recurrence relation here 
  differs only in the factor $2\A^2$, which does not affect the argument.
  This completes the proof.
  
  It is perhaps interesting that there is a direct calculation leading
  to the Pfaffian formula for the multidegree of $D_1$. Here we use a
  slightly different, analytic, language to emphasize the connection
  to matrix models.  We give the details of the calculation in the
  case $N$ even.

  The action of the torus $T$, and the moment map 
  $$ \Phi : \MNC \to \Lie(T)^*, \qquad 
  M\mapsto \pi \sum_{i,j} \left| M_{ij}\right|^2 (\A+z_i-z_j) $$
  both restrict to $D_1$. 
  Using the matrix $Z={\rm diag}(z_1,\ldots,z_N)$, we can rewrite this
  as $\Phi(M) = \pi (\A\Tr MM^\dagger+\Tr Z[M,M^\dagger])$.

  Writing $c$ for $\prod_{i, j} (\A+z_i-z_j)$,
  a formal application of the push-pull formula leads to the formula
  $$
  \mdeg D_1 = c \int_{M\in D_1} d\mu(M) 
    \e{-\pi (\A\Tr MM^\dagger+\Tr Z[M,M^\dagger])} 
  $$
  where the measure $d\mu(M)$ on $D_1$ %and $d\mu_0(M)$ on $\MNC$ are
  is derived from the flat metric $\sum_{i,j}|M_{ij}|^2$.
  %  The denominator equals $\prod_{i,j} (\A+z_i-z_j)^{-1}\equiv c^{-1}$.
  It is not our intention to provide a rigorous justification of the
  above, but we will show that it leads to the correct formula we have
  already justified by other means.

  There is a decomposition of $M\in E_N$ as $M=\Omega M' \Omega^\dagger$
  where $\Omega$ is unitary and $M'$ has the $n\times n$ block structure
  $$M'=\begin{pmatrix}0&0\\ X&0\end{pmatrix}$$
  and $X$ is a diagonal matrix: $X={\rm diag}(x_1,\ldots,x_n)$ with $x_i\ge 0$.
  To find such a decomposition, first obtain the obvious block decomposition
  with $X$ arbitrary  ($\Img M \subset \Ker M$,
  $\dim \Ker M\ge n$), then use the standard fact that for any
  $n\times n$ complex matrix $X$ there exist $n\times n$
  unitary matrices $V$, $W$ such that $VXW^\dagger$ is diagonal positive.
  
  Noting that $\Tr MM^\dagger=\sum_{i=1}^n x_i^2$ we perform the change
  of variables in the integral. The measure in the new variables must be
  carefully computed by setting $\Omega=1+i d\Omega$ with
  $d\Omega=\begin{pmatrix}H_{11}&H_{12}\\ H_{21}&H_{22}\end{pmatrix}$
  Hermitian, $x'_i=x_i+dx_i$, expanding the metric $\sum_{i,j} |M_{ij}|^2$ 
  at first order in $d\Omega$ and $d x_i$, and finally taking the square
  root of its determinant. The diagonal parts $H_{11}$, $H_{22}$
  contribute the usual factors $\prod_{i=1}^n x_i \prod_{i<j} (x_i^2-x_j^2)^2$,
  but remarkably the part $H_{12}$ contributes
  $\prod_{i=1}^n x_i^2 \prod_{i<j} (x_i^2+x_j^2)^2$, so that this
  recombines into
  $$
  \mdeg D_1 = c \int {d \Omega\over (2\pi)^n n!} \prod_{i=1}^n d x_i
  x_i^3 \e{-\pi \A x_i^2} \Delta^2(x_i^4) \e{-\pi\Tr
    Z\Omega(XX^\dagger-X^\dagger X)\Omega^\dagger}
  $$
  where the factor $(2\pi)^n n!$ comes from the non-uniqueness of
  the decomposition, and $\Delta(\cdot)$ is the Vandermonde
  determinant: $\Delta(x_i^4)=\prod_{i<j} (x_i^4-x_j^4)$.  The
  integral over the unitary group is the Harish
  Chandra--Itzykson--Zuber integral \cite{HC,IZ} (see also
  \cite{ZJZ}). The diagonal matrices $Z$ and $XX^\dagger-X^\dagger X$ have 
  entries respectively $z_j$, $j=1,\ldots,N$, and $\pm x_i^2$, $i=1,\ldots,n$; 
  we write the latter as $XX^\dagger-X^\dagger X={\rm diag}(\epsilon x_i^2)$, 
  $(i,\epsilon)\in\{1,\ldots,n\}\times\{-1,+1\}$. We thus find
  $$
  \mdeg D_1= %{2^n \pi^{2n^2} \over n!\prod_{k=0}^{2n-1} k!} 
  c\, {(2\pi)^N\over (2\pi)^n n!}
  \int_0^\infty \prod_{i=1}^n d x_i x_i^3 \e{-\pi \A x_i^2} \Delta^2(x_i^4)
  {\det\left(\e{\epsilon \pi z_j x_i^2}\right)\over 
    \Delta(z_j)\Delta(\epsilon x_i^2)}
  $$
  $\Delta(\epsilon x_i^2)
  = \prod_{i<j} (x_i^2-x_j^2)^2(x_i^2+x_j^2)^2\prod_i(2x_i^2)$, 
  so that one can simplify and compute
  \begin{eqnarray*}
    \mdeg D_1&=& 
    c\, {\pi^n\over n!}
    \int_0^\infty \prod_{i=1}^n d x_i x_i \e{-\pi \A x_i^2}
    {\det(\e{\pi \epsilon z_j x_i^2})\over \Delta(z_j)}
    \\
    &=&c\, {\pi^n\over n!} \sum_{\sigma\in S_N} (-1)^\sigma
    {1\over\Delta(z_j)}
    \int_0^\infty \prod_{i=1}^n d x_i x_i 
    \e{-\pi x_i^2(\A+z_{\sigma(2i-1)}-z_{\sigma(2i)})}    \\
    &=& \A^N \prod_{i<j} {(\A-z_i+z_j)(\A-z_j+z_i)\over (z_i-z_j)} 
    {\rm Pf} \left( {1\over \A+z_j-z_i}-{1\over \A+z_i-z_j}\right)
  \end{eqnarray*}
  Since 
  ${1\over \A+z_j-z_i}-{1\over \A+z_i-z_j}={2(z_i-z_j)\over \A-(z_i-z_j)^2}$, 
  we obtain the desired expression; it differs from that of corollary
  \ref{cor:totalmdeg} by a factor of $2^n \A^N$. The power of $\A$ is simply 
  due to the different embedding space ($\MNC$ versus $\MNC_{\Delta=0}$).
  
  For $N$ odd the result of the computation of the integral is
  strictly identical; however to obtain the multidegree of $D_1$ one
  must take into account the multiplicity $2$, hence the factor $2^r$.
\end{proof}

\begin{Theorem}
  Let 
  $$ D_t := \{ M : (M_\leq + t M_>)^2 = 0 \}, \qquad t\neq 0 $$
  so each $D_t \iso D_1 = \{ M : M^2 = 0 \}$.  Define $D_0$ to be the
  flat limit $\lim_{t\to 0} D_t$.  Then the scheme $D_0$ is supported
  on $\cup_\pi E_\pi$, and contains each $E_\pi$ with the same multiplicity
  $2^{n+r}$.
\end{Theorem}

\begin{proof}
  As explained in section \ref{ssec:algdegen}, the limit of the set of
  equations $(M_\leq + t M_>)^2 = 0$ as $t\to 0$ is the set $M\cp M = 0$. 
  However, these may not generate the limit ideal defining $D_0$. 
  So we can only infer a containment (of schemes), 
  $ D_0 \subseteq \{M\in \MNC : M\cp M = 0\}$. 
  While this latter scheme is bigger
  than $E$, it has the same support, so as sets $D_0 \subseteq E$.

  Since $D_1$ is irreducible and hence equidimensional, the flat limit
  $D_0$ is also equidimensional, so it is supported on $E$'s components
  of top dimension, $\cup_\pi E_\pi$. (Remember that we conjecture that
  $E$ has no other components, but even if it does they're not in $D_0$.)

  Consequently 
  $$ \mdeg D_1 = \mdeg D_0 = \sum_\pi c_\pi\ \mdeg E_\pi $$
  for some coefficients $\{c_\pi \in \naturals\}$, where $c_\pi$ is
  the multiplicity of $E_\pi$ in $D_0$.

  However, we already know $\mdeg D_1$ from proposition \ref{prop:degD1}
  and $\sum_\pi \mdeg E_\pi $ from corollary \ref{cor:totalmdeg}, from which 
  we see that taking $c_\pi \equiv 2^{n+r}$ gives a solution. 
  To know it's the right one, it is enough to show that the polynomials 
  $\{ \mdeg E_\pi\}$ are linearly independent over $\integers$.

  Let $\sum_\pi d_\pi\ \mdeg E_\pi = 0$ be a linear relation among them.
  By theorem \ref{thm:main}, we also know $\sum_\pi d_\pi\ \Psi_\pi = 0$.
  Let $\rho$ be a link pattern.
  By \cite[lemma 2]{DFZJ}, the specialization of $\Psi_\pi$ at $\A=0$, 
  $z_i = z_{\rho(i)}, i=1\ldots N$ is nonzero if and only if $\pi=\rho$,
  allowing us to pick out the $d_\rho$ term and show $d_\rho=0$.
  
  Hence the $\{\mdeg E_\pi\}$ are linearly independent, and the
  multiplicities are all $2^{n+r}$.
\end{proof}

The power of $2$ in the multiplicities on $D_0$ can be loosely blamed
on the ``missing'' equations $M_{ii}=0$ included in the definition of $E$.

We conjecture that $D_0 = \{M\in \MNC : M\cp M = 0\}$ as schemes,
which would imply our earlier conjecture that $E$ is equidimensional, 
in that $D_0$ is the flat limit of a variety and hence equidimensional.

\section{An additional circle action}\seclabel{extra}

Throughout this paper we claimed to be working with the action of an
$N+1$-dimensional torus $T$ on the scheme $E$. Since the
$1$-dimensional subtorus of $T$ consisting of scalar matrices acts
trivially, it is really more honest to consider this an action of
the $N$-dimensional quotient torus. The corresponding statement for
the multidegrees is that while we considered our multidegrees as polynomials
in $\A,z_1,\ldots,z_N$, they can all be written as polynomials in the $N$ 
expressions $\A,z_1-z_2,\ldots,z_{N-1}-z_N$. (Of course, for the multidegrees
of subschemes of $\MNCdzero$, we have an even better statement -- they
are polynomials with positive coefficients in the weights 
$\{\A+z_i-z_j : i\neq j \}$ of $\MNCdzero$.)

Abstractly, we should expect that $(\MNC,\cp)$ has an extra degree of symmetry
beyond that of $(\MNC,\times)$, in that $\cp$ is the multiplication on the 
degenerate fiber of a $1$-parameter family (see %the $t\cdot M$ action in 
section \ref{ssec:algdegen}). It is easy to write down this bigger action: 
define
$$ (\alpha,\omega_1,\ldots,\omega_N) \cdot e^{ik} 
     := \alpha \left( \prod_{j: \cyc{i\leq j < k}} \omega_j\right) e^{ik}. $$
If $\omega_i = \zeta_i \zeta_{i+1}^{-1}$
for each $i$, then the action of $(\alpha,\omega_1,\ldots,\omega_N)$ is just
conjugation by the diagonal matrix $\diag(\zeta_1,\ldots,\zeta_N)$, followed by
rescaling by $\alpha$. Hence this extends the $T$-action.

It also is easy to check that the action of the subgroup with $\alpha=1$
preserves the product $\cp$ on $\MNC$. First,
\begin{eqnarray*}
 (1,\omega_1,\ldots,\omega_N) \cdot (e^{hj} \cp e^{km})  
&=& (1,\omega_1,\ldots,\omega_N) \cdot (\delta_{jk} [\cyc{h\leq j\leq m}] e^{hm})
\\
&=& \left( \prod_{i: \cyc{h\leq i< m}} \omega_i\right) 
\delta_{jk} [\cyc{h\leq j\leq m}]e^{hm}.
\end{eqnarray*}
Now notice that $\cyc{h\leq j\leq m}$ implies that
$$ \prod_{i: \cyc{h\leq i< m}} \omega_i 
= \prod_{i: \cyc{h\leq i< j}} \omega_i   \prod_{l: \cyc{j\leq l< m}} \omega_l $$
which is what we need to establish
$$  (1,\omega_1,\ldots,\omega_N) \cdot (e^{hj} \cp e^{km})  
= \left( (1,\omega_1,\ldots,\omega_N) \cdot e^{hj} \right) \cp 
\left((1,\omega_1,\ldots,\omega_N) \cdot e^{km}  \right) $$
when both sides are nonzero.
Since this action with $\alpha=1$ preserves $\cp$, it preserves the scheme $E$,
and $\alpha$ is just acting by rescaling $E$. 

\newcommand\barT{{\overline T}}
\newcommand\barPsi{{\overline \Psi}}
\newcommand\wtot{w_{tot}}
Call this bigger torus $\barT$, and use $(\A,w_1,\ldots,w_N)$ for
the obvious basis of its weight lattice. Then we get the following equation
on the $\barT$-- and $T$--multidegrees of an affine scheme
$X \subseteq \MNCdzero$:
$$ \barT\!-\!\mdeg\ X |_{w_i = z_i - z_{i+1}} = T\!-\!\mdeg\ X $$
where each $w_i$ has been specialized to $z_i - z_{i+1}$.
The kernel of this specialization is generated by 
$$ \wtot := \sum_{i=1}^N w_i. $$

%The proof of proposition \ref{prop:recurrence}, 
%which applied equally well to $\barT$-multidegrees, actually proves 
%\begin{equation}\label{eqn:transpb}
%  \mdeg E_{\pi} + \mdeg E_{f_i\cdot\pi}
%  =-\ \frac{2a - b_i + \btot}{a - b_i + \btot}\ 
%  \der_i\left( (a - b_i + \btot)\ \mdeg E_{\pi}\right).
%\end{equation}
%because the $\barT$-weights on the $(i,i+1)$ entries of $M$ and $M^2$
%are $a-b_i+\btot$ and $2a-b_i+\btot$. The base case,
%$\pi_0(i)=i+n \bmod 2n$ for $i\leq 2n$, and $\pi_0(N)=N$ if N is odd,
%is also just as easy as in proposition \ref{prop:pizero}:
%\begin{equation*} \label{eqn:pizerob}
%  \mdeg E_{\pi_0}
%  = \prod_{i=1\ldots N \atop j:\ \cyc{i<j<i+n}} 
%  \bigg(a+\sum_{k:\ \cyc{i\leq k<j}} b_k\bigg)
%  \left( \prod_{i=n+1}^N  \bigg(a+\sum_{k:\ \cyc{i\leq k<i+n}} b_k\bigg) 
%  \right)^r
%\end{equation*}
%Lacking a geometric proof of corollary \ref{cor:smallarch}, we don't
%know how to lift that formula to one about $\barT$-multidegrees.

The generalization of theorem \ref{thm:main} will be discussed in 
\cite{DFKZJ}.


\begin{thebibliography}{10}

\bibitem[BS]{BS} R. Bott, H. Samelson,
  The cohomology ring of $G/T$.
  Proc. Nat. Acad. Sci. U. S. A. 41 (1955), 490--493.

\bibitem[Br]{Br}  M. Brion,
  Equivariant cohomology and equivariant intersection theory.
  Notes by Alvaro Rittatore. NATO Adv. Sci. Inst. Ser. C Math. Phys. Sci., 514,
  Representation theories and algebraic geometry (Montreal, PQ, 1997), 1--37.
  {\tt math.AG/9802063}

\bibitem[dGN]{dGN} J. de Gier, B. Nienhuis,
  Brauer loops and the commuting variety. 
  J. Stat. Mech. (2005) P01006.
  {\tt math.AG/0410392}
  
\bibitem[DFZJ]{DFZJ} P. Di Francesco, P. Zinn-Justin,
  Inhomogeneous model of crossing loops and multidegrees of some
  algebraic varieties. 
  {\tt math-ph/0412031}

\bibitem[DFKZJ]{DFKZJ} P. Di Francesco, A. Knutson, P. Zinn-Justin,
  The Brauer loop scheme and orbital varieties,
  in preparation.

\bibitem[Fu]{Fu} W. Fulton, 
  Flags, Schubert polynomials, degeneracy loci,
  and determinantal formulas. Duke Math. J. \textbf{65} (1992), no.~3,
  381--420.

\bibitem[GS]{GS} V. Guillemin, S. Sternberg,
  {Supersymmetry and equivariant de Rham theory.}
  Springer-Verlag, Berlin, 1999. 

\bibitem[HC]{HC} Harish Chandra, 
  Differential operators on a semisimple Lie algebra.
  Amer. J. Math. {\bf 79} (1957) 87--120.

\bibitem[IZ]{IZ} C. Itzykson and J.-B. Zuber, 
  The planar approximation. II.
  J. Math. Phys. {\bf 21} (1980) 411--421.

\bibitem[Jo]{Jo} A. Joseph, 
  {On the variety of a highest weight module}.
  J. Algebra \textbf{88} (1984), no.~1, 238--278.
  
\bibitem[KR]{KR} V. G. Kac and A. K. Raina, 
  {Bombay lectures on highest weight representations of infinite
    dimensional Lie algebras,} (Lecture 9). Advanced Series in
  Mathematical Physics Vol. 2, World Scientific.

\bibitem[Kn]{Kn}   A.~Knutson,
  Some schemes related to the commuting variety. 
  To appear in the Journal of Algebraic Geometry.
  {\tt math.AG/0306275}

\bibitem[KMY]{KMY} A. Knutson, E. Miller, A. Yong,
  Gr\"obner geometry of vertex decompositions and of flagged tableaux.
  preprint 2005.
  {\tt math.AG/0502144}

\bibitem[M]{M}  A. Melnikov,
  $B$-Orbits in Solutions to the Equation $X^2 = 0$ in Triangular Matrices.
  Journal of Algebra 223, 101--108 (2000).

\bibitem[MR]{MR} M. J.~Martins and P. B.~Ramos, 
  The Algebraic Bethe Ansatz for rational braid-monoid lattice models.
  Nucl. Phys. B500 (1997) 579--620.
  {\tt hep-th/9703023}

\bibitem[MNR]{MNR} M. J.~Martins, B.~Nienhuis and R.~Rietman, 
  An Intersecting Loop Model as a Solvable Super Spin Chain.
  Phys. Rev. Lett. 81 (1998) 504--507.
  {\tt cond-mat/9709051}

\bibitem[MS]{MS}   E.~Miller and B.~Sturmfels, 
  {Combinatorial commutative algebra}.
  Graduate Texts in Mathematics, vol.~227, Springer--Verlag, New York, 2004.

\bibitem[NR]{NR} B. Nienhuis and R. Rietman, 
  A solvable loop model with intersections, preprint IFTA-92-35;
  R. Rietman, Yang--Baxter equations, Hyperlatices and a loop model,
  PhD thesis (unpublished)

\bibitem[Ro]{Ro}   W. Rossmann,
  Equivariant multiplicities on complex varieties.
  Orbites unipotentes et repr\'esentations, III.
  Ast\'erisque No. 173-174, (1989), 11, 313--330.

\bibitem[St]{St}  E. Strickland,
  On the variety of projectors.
  J. Algebra 106 (1987), no. 1, 135--147.

\bibitem[We]{We}  J. Weyman,
  Two results on equations of nilpotent orbits. 
  J. Algebraic Geom. 11 (2002), no. 4, 791--800.
  {\tt math.AG/0006232}

\bibitem[ZJZ]{ZJZ}  P. Zinn-Justin, J.-B. Zuber,
  On some integrals over the $U(N)$ unitary group and their large $N$ limit.
  J.Phys. A36 (2003) 3173--3194.
  {\tt math-ph/0209019}

\end{thebibliography}
\end{document}